\documentclass{amsart}
\usepackage{amssymb,citesort,texdraw}

\theoremstyle{plain}
   \newtheorem{thm}{Theorem}[section]
   \newtheorem{prop}[thm]{Proposition}
   \newtheorem{lem}[thm]{Lemma}
   
\theoremstyle{definition}
   \newtheorem{df}[thm]{Definition}
\theoremstyle{remark}

\renewcommand{\theenumi}{(\alph{enumi})}
\renewcommand{\labelenumi}{\theenumi}
\numberwithin{equation}{section}
\everytexdraw{\drawdim mm \setunitscale 1 \linewd 0.3 \move(0 0)}
\newcommand{\rcir}{\bsegment \move(1 0) \lcir r:1 \savepos(2 0)(*px *py)
                   \esegment \move(*px *py)}

\newcommand{\Z}{\mathbf{Z}}
\newcommand{\ZP}{\mathfrak{Z}}
\newcommand{\C}{\mathbf{C}}
\newcommand{\R}{\mathbf{R}}
\newcommand{\Hf}{H}
\newcommand{\Ct}{C}
\newcommand{\I}{I}
\newcommand{\m}{\underline{m}}
\newcommand{\g}{\mathfrak{g}}
\newcommand{\h}{\mathfrak{h}}
\newcommand{\ali}{\alpha_i}
\newcommand{\alj}{\alpha_j}
\newcommand{\eik}{e_{i,k}}

\newcommand{\eil}{e_{i,l}}
\newcommand{\ejl}{e_{j,l}}
\newcommand{\fik}{f_{i,k}}

\newcommand{\fjl}{f_{j,l}}
\newcommand{\cntr}{\mathfrak{z}}
\newcommand{\defi}{\emph}
\newcommand{\brbin}[2]{\genfrac{\{}{\}}{0pt}{}{#1}{#2}}
\DeclareMathOperator{\ad}{ad}
\DeclareMathOperator{\id}{id}
\DeclareMathOperator{\im}{Im}
\DeclareMathOperator{\adt}{\widetilde{ad}}
\DeclareMathOperator{\str}{str}
\DeclareMathOperator{\tr}{tr}

\begin{document}

\title[CENTER AND R-MATRIX FOR QUANTIZED BORCHERDS SUPERALGEBRAS]
      {Center and Universal R-matrix\\
       for\\
       Quantized Borcherds Superalgebras}
\author{Jin Hong}
\address{Department of Mathematics\\
         Seoul National University\\
         Seoul 151-742, Korea}
\thanks{Supported in part by GARC-KOSEF at Seoul National University}
\email{jhong@math.snu.ac.kr}
\begin{abstract}
We construct a nondegenerate symmetric bilinear form on quantized
enveloping algebras associated to Borcherds superalgebras.
With this, we study its center and its universal R-matrix.
\end{abstract}
\maketitle

\section{Introduction}

Quantized enveloping algebras for Kac-Moody algebras were introduced
independently by Drinfel'd(\cite{Drf}) and Jimbo(\cite{Jim})
in their studies of the Yang-Baxter equation.
The Kac-Moody algebras were generalized(\cite{Bch}) by Borcherds to
accommodate for his study of the monstrous moonshine(\cite{mm}).
And quantized version of the enveloping algebras for
Borcherds algebras(\cite{Kg95}) was soon studied.
There are also superalgebra versions of these algebras(\cite{Kang96}).
We shall study the structure of the center and find the R-matrix
for the quantized Borcherds superalgebras.

Much work has been done on the center of quantized enveloping algebras
for finite dimensional semisimple
Lie algebras(\cite{Bau,Dr90,JL,Rsh,RTF,Ro,Ta91}),
and there are Kac-Moody(\cite{Pa}) and Borcherds(\cite{KgTa}) versions
also.
We will mainly follow \cite{Ta91} and \cite{KgTa} to find the
center for quantized Borcherds superalgebras.

As for the universal R-matrix, the quantum double construction by
Drinfel'd(\cite{Dr87}) gives its existence for any Hopf algebra satisfying
some conditions.
Even though there is a quantum double construction for $\Z_2$-graded
Hopf algebras(\cite{Gould}), we do not use it in this paper.
Instead, we explicitly construct a universal R-matrix and show that it
satisfies the Yang-Baxter equation.

The paper is organized as follows.
In Section~\ref{QDBS}, we define the quantized Borcherds superalgebras
and give it a Hopf algebra structure. The triangular decomposition will
also be mentioned.
In Section~\ref{RU}, the character formula for highest weight
representations will
be given and we prove a lemma that will be used in later sections.
The next section is devoted to providing the quantized Borcherds
superalgebras with a bilinear form and proving its nondegeneracy.
In Section~\ref{HC}, we define the Harish-Chandra homomorphism, show its
injectivity, and prove some properties concerning its image.
Information on the center of the quantized Borcherds superalgebra will
be obtained in Section~\ref{CU}.
The last section will give the universal R-matrix and show that it
satisfies the Yang-Baxter equation.

\section{Quantum Deformation of Borcherds Superalgebras}\label{QDBS}

In this section, we define the quantized Borcherds superalgebras and
give it a Hopf algebra structure.

Let $\I$ be a countable index set. A matrix $A=(a_{i,j})_{i,j\in \I}$
with entries in the real numbers is a \defi{Borcherds-Cartan} matrix if
\begin{itemize}
\item $a_{i,i} = 2$ or $a_{i,i} \leq 0$ for all $i\in\I$,
\item $a_{i,j} \leq 0$ if $i \neq j$ and
      $a_{i,j} \in \Z$ if $a_{i,i} = 2$,\label{hihi}
\item $a_{i,j} = 0$ if and only if $a_{j,i} = 0$.
\end{itemize}
If there exists a diagonal matrix $D = \text{diag}(s_i|i\in\I, s_i > 0)$
such that $DA$ is symmetric, then $A$ is said to be \defi{symmetrizable}.
If a symmetrizable Borcherds-Cartan matrix $A$ further satisfies the
constraints,
\begin{itemize}
\item $a_{i,j} \in \Z$,
\item $a_{i,i} \in 2\Z$,
\item $s_i \in \Z_{>0}$,
\end{itemize}
for all $i,j \in \I$, then it is said to be \defi{integral}.

A complex matrix $C = (\theta_{i,j})_{i,j\in\I}$ is a \defi{coloring matrix}
if $\theta_{i,j}\theta_{j,i} = 1$ for all $i,j \in\I$. Necessarily,
$\theta_{i,i} = \pm 1$ and we say $i$ is \defi{even} when $\theta_{i,i} = 1$,
\defi{odd} when $\theta_{i,i} = -1$.
A Borcherds-Cartan matrix $A$ is \defi{colored by $C$}
if for every $i\in\I$ such that $a_{i,i}=2$ and $\theta_{i,i} = -1$
we have $a_{i,j}\in 2\Z$ for all
$j\in\I$.

Throughout this paper, we shall assume that $A$ is a symmetrizable integral
Borcherds-Cartan matrix which is colored by a coloring
matrix $C$.

Let $\I^{\text{re}} = \{ i\in\I|a_{i,i}=2 \}$ and
$\I^{\text{im}} = \{ i\in\I|a_{i,i}\leq 0 \}$.
Also let $\m = (m_i|i\in\I)$ be a collection of positive integers such that
$m_i=1$ for all $i\in\I^{\text{re}}$. We call $\m$ the \defi{charge} of the
Borcherds-Cartan matrix $A$.

For a symmetrizable integral Borcherds-Cartan matrix $A$, which is
colored by a coloring matrix $C$, we denote by $\g(A,\m,C)$ the
Borcherds superalgebra of charge \m. (See~\cite{Kang96}.)

We set $P^{\vee} = (\bigoplus_{i\in\I}\Z h_i)\oplus(\bigoplus_{i\in\I}\Z d_i)$
and let $\h = \C\otimes_{\Z}P^{\vee}$ be the complex vector space with
basis $\{h_i,d_i|i\in\I\}$. For $i\in\I$, we define $\ali$ in the dual
space $\h^{*}$ of $\h$ by setting $\ali(h_j) = a_{j,i}$ and
$\ali(d_j) = \delta_{i,j}$.
Since $A$ is assumed to be symmetrizable, there exists a nondegenerate
symmetric bilinear form $(\ |\ )$ on $\h$ given by
$(s_ih_i|h) = \ali(h)$ and $(d_i|d_j) = 0$ for $i,j\in\I$, $h\in\h$.

The free abelian group $Q = \bigoplus_{i\in\I}\Z\ali$ generated by the
$\ali$ $(i\in\I)$ is called the \defi{root lattice} associated to $A$.
Let $Q^+ = \sum_{i\in\I}\Z_{\geq 0}\ali$ and $Q^- = -Q^+$.
The coloring matrix $C=(\theta_{i,j})$ gives
rise to a complex valued mapping
$\theta : Q \times Q \longrightarrow \C^{\times}$ satisfying
\begin{itemize}
\item $\theta(\ali,\alj) = \theta_{i,j}$,
\item $\theta(\alpha,\beta+\gamma)
        = \theta(\alpha,\beta)\theta(\alpha,\gamma)$,
\item $\theta(\alpha+\beta,\gamma)
        = \theta(\alpha,\gamma)\theta(\beta,\gamma)$,
\end{itemize}
for all $\alpha, \beta, \gamma \in Q$.

We define the binomial coefficients by:\\
$\{n\}_{q_i} =%
\frac{\theta_{i,i}^{n}q_i^n - q_i^{-n}}{\theta_{i,i}q_i - q_i^{-1}}$,
$\{n\}_{q_i}! = \overset{n}{\underset{t=1}{\prod}}\{t\}_{q_i}$, and
$\brbin{m}{n}_{q_i} = %
\frac{\{m\}_{q_i}!}{\{n\}_{q_i}!\{m-n\}_{q_i}!}$,\\
where $\{0\}_{q_i}!=1$ and $q_i = q^{s_i}$.

We let $\xi_i = q_i - q_i^{-1}$ and $K_i = q^{s_i h_i}$.

\begin{df}[\cite{Kang96}]\label{qtsupalg}
Suppose $\g = \g(A,\m,C)$ is the Borcherds superalgebra of charge \m\
determined by the symmetrizable integral Borcherds-Cartan matrix $A$ which
is colored by a coloring matrix $C$. Let $q$ be an
indeterminate. Then the \defi{quantized Borcherds superalgebra} $U_q(\g)$
associated to $\g$ is the
associative algebra over $\C(q)$ with 1, generated by the elements $q^h$
$(h\in P^{\vee})$, $\eik$, $\fik$ $(i\in\I, k=1,2,\cdots,m_i)$
with the defining relations:
{ 
\renewcommand{\theenumi}{(R\arabic{enumi})}
\renewcommand{\labelenumi}{\theenumi}
\begin{enumerate}
\item $q^0 = 1$, $q^h q^{h'} = q^{h+h'}$\quad for $h,h' \in P^{\vee}$,
\item $q^h \eik q^{-h} = q^{\ali(h)}\eik$\quad
      for $h\in P^{\vee}$, $i\in\I$, $k=1,2,\cdots,m_i$,
\item $q^h \fik q^{-h} = q^{-\ali(h)}\fik$\quad
      for $h\in P^{\vee}$, $i\in\I$, $k=1,2,\cdots,m_i$,
\item $\eik\fjl - \theta_{j,i}\fjl\eik = 
      \delta_{i,j}\delta_{k,l}\frac{1}{\xi_i}(K_i-K_i^{-1})$,\\[2mm]
      for $i,j\in\I$, $k=1,2,\cdots,m_i$, $l=1,2,\cdots,m_j$,
\item $\overset{1-a_{i,j}}{\underset{n=0}{\sum}}%
      (-1)^n \theta_{i,j}^n\theta_{i,i}^{n(n-1)/2}%
      \brbin{1-a_{i,j}}{n}_{q_i}\eik^{1-a_{i,j}-n}\ejl\eik^n
      = 0$\\[2mm]
      if $a_{i,i}=2$ and $i \neq j$,
\item $\overset{1-a_{i,j}}{\underset{n=0}{\sum}}%
      (-1)^n \theta_{i,j}^n\theta_{i,i}^{n(n-1)/2}%
      \brbin{1-a_{i,j}}{n}_{q_i}\fik^{1-a_{i,j}-n}\fjl\fik^n
      = 0$\\[2mm]
      if $a_{i,i}=2$ and $i \neq j$,
\item $\eik\ejl - \theta_{i,j}\ejl\eik = 0$\quad if $a_{i,j} = 0$,
\item $\fik\fjl - \theta_{i,j}\fjl\fik = 0$\quad if $a_{i,j} = 0$.
\end{enumerate}
} 
\end{df}

\begin{prop}[\cite{Kang96}]
The algebra $U_q(\g)$ has a Hopf algebra structure
with comultiplication $\Delta$,
counit $\varepsilon$, and antipode $S$ defined by:
\begin{align}
\Delta(q^h) &= q^h \otimes q^h,\\
\Delta(\eik) &= \eik \otimes 1 + K_i\otimes \eik,\\
\Delta(\fik) &= \fik \otimes K_i^{-1} + 1 \otimes \fik,
\end{align}
\begin{align}
\varepsilon(q^h) &= 1,\\
\varepsilon(\eik) &= 0,\\
\varepsilon(\fik) &= 0,
\end{align}
\begin{align}
S(q^h) &= q^{-h},\\
S(\eik) &= -K_i^{-1}\eik,\\
S(\fik) &= -\fik K_i,
\end{align}
for $h\in P^{\vee}$, $i\in\I$, $k=1,2,\cdots,m_i$.
\end{prop}

We denote by $U^0$ the subalgebra of $U = U_q(\g)$ generated by $q^h$
for $h\in P^{\vee}$ and $U^+$ (respectively, $U^-$) the subalgebra of $U$
generated by the elements $\eik$ (respectively, $\fik$) for $i\in\I$,
$k = 1,2,\cdots,m_i$. We also denote by $U^{\geq 0}$ (respectively,
$U^{\leq 0}$) the subalgebra of $U$ generated by the elements $q^h$ and $\eik$
(respectively, $\fik$) for $h\in P^{\vee}$, $i\in\I$, $k=1,2,\cdots,m_i$.
For each $\beta\in Q$, let
\begin{equation}
U_{\beta} = \{ x\in U \mid %
q^h x q^{-h} = q^{\beta(h)}x \text{ for all } h\in P^{\vee} \}.
\end{equation}
We similarly define $U_{\pm\beta}^{\pm}$, $U_{\pm\beta}^{\geq0}$, and 
$U_{\pm\beta}^{\leq0}$ for $\beta\in Q^+$.
We then have:
\begin{prop}[\cite{Kang96}]\label{tridcmp}\hfill
\begin{enumerate}
\item $U \cong U^-\otimes U^0 \otimes U^+$.
\item $U^0 = \bigoplus_{h\in P^{\vee}} \C q^h$.
\item $U^{\pm} = \bigoplus_{\beta\in Q^+} U_{\pm\beta}^{\pm}$.
\item \textnormal{(R5)} and \textnormal{(R7)} \textup{(}respectively,
\textnormal{(R6)} and \textnormal{(R8)}\textup{)} are the fundamental
relations for $U^+$ \textup{(}respectively, $U^-$\textup{)}.
\end{enumerate}
\end{prop}

We give $Q^+$ a partial ordering by setting $\lambda \geq \mu$ if
and only if $\lambda - \mu \in Q^+$.
We will also use the notation $K_\gamma = \prod K_i^{n_i}$\label{nata} for
$\gamma = \sum n_i \ali \in Q$.

\section{Representations of $U_q(\g)$}\label{RU}

For $i\in\I$ define the $\C$-linear functionals $\Lambda_i \in \h^*$
by:
\begin{equation}
\Lambda_i(h_j) = \delta_{i,j},\qquad\Lambda_i(d_j) = 0,
\qquad\text{for all } j\in\I.
\end{equation}
Define the lattices:
\begin{align}
P &= \{ \lambda \in \h^* |%
          \lambda(h_i), \lambda(d_i) \in \Z, \forall i\in\I \},\\
\overline{P} &= %
          (\bigoplus_{i\in\I}\Z\ali)\oplus(\bigoplus_{i\in\I}\Z\Lambda_i).
\end{align}
$P$ is called the \defi{weight lattice} of $\g$.
An element $\lambda\in P$ is said to be a \defi{dominant integral weight} if
\begin{align}
\lambda(h_i) &\in \Z_{\geq0} \text{ for all } i\in\I^{\text{re}},\\
\lambda(h_i) &\in 2\Z_{\geq0}
		\text{ for all } i\in\I^{\text{re}} \cap\I^{\text{odd}},
\end{align}
where $\I^{\text{odd}}$ denotes the set of $i\in\I$ such that
$\theta_{i,i}=-1$. Let $P^+$ denote the set of all dominant integral weights.

Set $\bar{\h}^* = \C\otimes_\Z \overline{P}$.
Then the nondegenerate symmetric bilinear form on $\h$ gives an isomorphism
between $\h$ and $\bar{\h}^*$ hence also induces a bilinear form
on $\bar{\h}^*$.
We may extend this bilinear form to a symmetric bilinear form on $\h^*$.
We extend it so that it satisfies $(\lambda|\ali) = \lambda(s_i h_i)$ and
$(\lambda|\Lambda_i) = \lambda(s_i d_i)$ for every $\lambda\in\h^*$.
Write $\lambda\bot\mu$ if $(\lambda|\mu)=0$.

For each $i\in\I$ such that $a_{i,i}\neq0$,
we define the \defi{simple reflection}
$r_i \in \text{GL}(\h^*)$ on $\h^*$ by
\begin{equation}
r_i(\lambda) = \lambda - \frac{2}{a_{i,i}}\lambda(h_i)\ali.
\end{equation}
The subgroup $W$ of $\text{GL}(\h^*)$ generated by
$r_i$ ($i\in\I^{\text{re}}$)
is called the \defi{Weyl group} of $\g(A,\m,C)$. We denote by
$l : W \longrightarrow \Z_{\geq0}$ the natural length function.

Let $R$ be the family of all imaginary simple roots,
each root occurring as many times as its multiplicity,
i.e. $m_i$ times for $\ali$.
For $\lambda\in P^+$, define $R(\lambda)$ to be the set of all
$\mu = \sum^r_{j=1}\alpha_{i_j} + \sum^s_{k=1}l_{i_k}\beta_{i_k} \in Q^+$,
where $\alpha_{i_j}$ (resp. $\beta_{i_k}$) are distinct even
(resp. odd) roots in $R$, satisfying
\begin{itemize}
\item $\alpha_{i_j}\bot\lambda$, $\beta_{i_k}\bot\lambda$ for all $j,k$,
\item $\alpha_{i_j}\bot\beta_{i_k}$ for all $j,k$,
\item $\alpha_{i_j}\bot\alpha_{i_k}$, $\beta_{i_j}\bot\beta_{i_k}$
	for $j\neq k$,
\item $\beta_{i_k}\bot\beta_{i_k}$ if $l_{i_k}\geq 2$.
\end{itemize}
In particular, $0\in R(\lambda)$.
Suppose $\rho\in\h^*$ satisfies $\rho(h_i) = \frac{1}{2} a_{i,i}$ for all
$i\in\I$.

\begin{prop}[\cite{Mi96,Ray}]\label{irhwt}
Let $\lambda \in P^+$.
Denote by $M^q(\lambda)$ the Verma module for $U_q(\g)$
with highest weight $\lambda$ and let $V^q(\lambda)$ be the irreducible
highest weight module over $U_q(\g)$ with highest weight $\lambda$.
Then,
\begin{equation}
\operatorname{ch}M^q(\lambda) = 
\frac{e^\lambda}{
\prod_{\alpha\in\Phi^-}
(1-\theta(\alpha,\alpha)e^\alpha)^{\theta(\alpha,\alpha)\dim\g_\alpha}
} = 
e^\lambda
\sum_{\beta\in Q^+} (\dim U^-_{-\beta})e^{-\beta},
\end{equation}
\begin{equation}
\operatorname{ch}V^q(\lambda) = 
\frac{
\sum_{w\in W, \mu \in R(\lambda)}
(-1)^{l(w)+\textnormal{ht}(\mu)} e^{w(\lambda+\rho-\mu)-\rho}
}{
\prod_{\alpha\in\Phi^-}
(1-\theta(\alpha,\alpha)e^\alpha)^{\theta(\alpha,\alpha)\dim\g_\alpha}
}.
\end{equation}
In this formula, $\Phi^-$ is the set of all negative roots.
\end{prop}

The following is a Corollary to this proposition.
\begin{lem}\label{UViso}
Let $\gamma = \sum_{i\in\I} n_i\ali \in Q^+$.
Suppose $\lambda\in P^+$, $\lambda(h_i) > 0$ for all
$i\in\I^{\textnormal{im}}$
and $\lambda(h_i) \geq n_i$ for all $i\in\I^{\textnormal{re}}$.
Then we have a linear isomorphism
$U^-_{-\gamma} \overset{\sim}{\longrightarrow} V^q(\lambda)_{\lambda-\gamma}$
given by $ u \longmapsto uv_\lambda$.
\end{lem}
\begin{proof}
$U^-_{-\gamma}\longrightarrow M^q(\lambda)_{\lambda-\gamma}$
is surjective, so
$U^-_{-\gamma} \longrightarrow V^q(\lambda)_{\lambda-\gamma}$
is also surjective.
Hence it suffices to show
$\dim U^-_{-\gamma} = \dim  V^q(\lambda)_{\lambda-\gamma}$.
Since $(\ali|\lambda) = \lambda(s_i h_i) = s_i\lambda(h_i) > 0$
for all $i\in\I^{\text{im}}$,
no nonempty subset $F$ of $R$ satisfies $F\bot\lambda$, and so
\begin{align}
\operatorname{ch}V^q(\lambda) &=
\frac{
\sum_{w\in W}
(-1)^{l(w)} e^{w(\lambda+\rho)-\rho}
}{
\prod_{\alpha\in\Phi^-}
(1-\theta(\alpha,\alpha)e^\alpha)^{\theta(\alpha,\alpha)\dim\g_\alpha}
}\\ &=
(\sum_{w\in W} (-1)^{l(w)}e^{w(\lambda+\rho)-\rho})
(\sum_{\beta\in Q^+}(\dim U^-_{-\beta})e^{-\beta})
\end{align}
Therefore, it suffices to show that if
$w(\lambda+\rho)-\rho -\beta = \lambda -\gamma$ for some $w\in W$,
$\beta\in Q^+$, then $w=1$.

We will show that if $w\neq1$, then
$\gamma+w(\lambda+\rho) -(\lambda+\rho) \not\in Q^+$ by using induction on the
length of $w$.\\
If $w=r_i$ $(i\in\I^{\text{re}})$, then
\begin{align*}
\gamma+r_i(\lambda+&\rho) - (\lambda + \rho)\\ 
&= \gamma + \lambda + \rho - (\lambda(h_i)+\rho(h_i))\ali - (\lambda + \rho)\\
&= \gamma - (\lambda(h_i) + 1)\ali \not\in Q^+.
\end{align*}
If $w=w'r_i$ $(i\in\I^{\text{re}})$ with $l(w)=l(w')+1$, then
\begin{align*}
\gamma+w(\lambda+&\rho) - (\lambda+\rho)\\
&= \gamma + w'r_i(\lambda+\rho) - (\lambda+\rho)\\
&= \gamma + w'(\lambda+\rho - (\lambda(h_i)+\rho(h_i))\ali)-(\lambda+\rho)\\
&= (\gamma+w'(\lambda+\rho)-(\lambda+\rho))-(\lambda(h_i)+1)w'\ali
\not\in Q^+.
\end{align*}
This completes the proof.
\end{proof}

\section{The Bilinear Form on $U_q(\g)$}

\subsection{The bilinear form on $U^{\geq0}\times U^{\leq0}$}
In this section, we define a bilinear form on $U^{\geq0}\times U^{\leq0}$
which is nondegenerate when restricted to
$U^+_\beta\times U^-_{-\beta}$, $\beta\in Q^+$.

For $\phi\in (U_\beta)^*$, $\psi\in (U_\gamma)^*$,
$x\in U_\beta$, and $y\in U_\gamma$, we define
$(\phi\otimes\psi)(x\otimes y) = \theta(-\gamma,\beta) \phi(x)\psi(y)$.
With this, and the Hopf algebra structure on $U_q(\g)$, we can give an
algebra structure to
$\bigoplus_{\alpha\in Q^+} (U^{\geq0}_{\phantom{-}\alpha})^*$ by setting
$(\phi_1\phi_2)(x) = (\phi_1\otimes\phi_2)(\Delta(x))$ for
$\phi_1,\phi_2\in \bigoplus_{\alpha\in Q^+} (U^{\geq0}_{\phantom{-}\alpha})^*$
and $x\in U^{\geq 0}$.
For $h\in P^\vee$ and $i\in\I$, $k=1,2,\cdots,m_i$, we define the linear
functionals $\phi_h,\psi_{i,k}\in%
\bigoplus_{\alpha\in Q^+} (U^{\geq0}_{\phantom{-}\alpha})^*$ by
\begin{alignat}{2}
\phi_h(xq^{h'}) &= \varepsilon(x)q^{-(h|h')} &\quad &(x\in U^+, h'\in P^\vee),\\
\psi_{i,k}(xq^h) &= 0
		&\quad &(x\in U^+_\beta, \beta\in Q^+\setminus\{\ali\}),\\
\psi_{i,k}(\eil q^h) &= \delta_{k,l}.
\end{alignat}

\begin{prop}\label{zetahom}
There exists an algebra homomorphism
\begin{equation}
\zeta : U^{\leq0} \longrightarrow
\bigoplus_{\alpha\in Q^+} (U^{\geq0}_{\phantom{-}\alpha})^*
\end{equation}
given by
\begin{alignat}{2}
\zeta(q^h) &= \phi_h & \quad&(h\in P^\vee),\\
\zeta(f_{i,k}) &= -\frac{1}{\xi_i}\psi_{i,k}&
\quad&(i\in\I, k=1,2,\cdots,m_i).
\end{alignat}
\end{prop}
\begin{proof}
By Proposition~\ref{tridcmp}, we have only to check that the relations
(R1), (R3), (R6), and (R8) are preserved under the map $\zeta$.
Other cases being easy, we just sketch the (R6) part.

Define $\eik^{(n)} = \eik^n/\{n\}_{q_i}!$.
We may check by induction on $n$ that
\begin{equation*}
\Delta(\eik^{(n)}) = 
\sum_{s+t=n} q_i^{st}\eik^{(s)}K_i^{t}\otimes\eik^{(t)}.
\end{equation*}
This show
\begin{equation*}
((\Delta\otimes1)\circ\Delta)(\eik^{(n)}) =
\sum_{r+s+t=n} q_i^{rs+st+tr}
\eik^{(r)}K_i^{s+t}\otimes \eik^{(s)}K_i^t\otimes \eik^{(t)}.
\end{equation*}
We again use induction to prove
\begin{equation*}
\psi_{i,k}^n(\eik^{(n)}) = (\theta_{i,i}q_i)^{\frac{n(n-1)}{2}}
\end{equation*}
With this, it is possible to show
\begin{equation*}
\psi_{i,k}^{N-n}\psi_{j,l}\psi_{i,k}^n (\eik^{(N-m)}\ejl \eik^{(m)}) =
\sum \theta_{i,i}^{A'}\theta_{i,j}^{B'}q_i^{C'}
\brbin{N-n}{\alpha}_{q_i}\brbin{n}{\beta}_{q_i},
\end{equation*}
with the summation over nonnegative integers $\alpha, \beta, \gamma, \delta$
such that $\alpha+\beta=N-m$, $\gamma+\delta=m$,
$\alpha+\gamma=N-n$, and $\beta+\delta=n$ and where
\begin{align*}
A' &= \beta\gamma+(N-n)n+\frac{1}{2}(N-n)(N-n-1)+\frac{1}{2}n(n-1),\\
B' &= \beta-\gamma+(N-n)-n,\\
C' &= \alpha\beta+\gamma\delta+2\beta\gamma+(\beta+\gamma)a_{i,j}
        +\frac{1}{2}(N-n)(N-n-1)+\frac{1}{2}n(n-1).
\end{align*}
Noting
\begin{equation*}
\brbin{N}{n}_{q_i}\brbin{N-n}{\alpha}_{q_i}\brbin{n}{N-m-\alpha}_{q_i} =
\brbin{N}{m}_{q_i}\brbin{N-m}{\beta}_{q_i}\brbin{m}{\delta}_{q_i},
\end{equation*}
we can calculate
\begin{align*}
\Big(\overset{1-a_{i,j}}{\underset{n=0}{\sum}}%
    (-1)^n \theta_{i,j}^n\theta_{i,i}^{n(n-1)/2}%
    \brbin{1-a_{i,j}}{n}_{q_i}\psi_{i,k}^{1-a_{i,j}-n}\psi_{j,l}\psi_{i,k}^n
\Big)
(\eik^{(N-m)}\ejl\eik^{(m)}) \\
\hfill=\sum_{n=0}^N \sum \theta_{i,i}^{A}\theta_{i,j}^{B}q_i^{C}
\brbin{N}{m}_{q_i}\brbin{N-m}{\beta}_{q_i}\brbin{m}{\delta}_{q_i}
\end{align*}
with the second summation over nonnegative integers
satisfying the same conditions as before and where
\begin{align*}
A &= \frac{1}{2}N(N-1)+m\beta+\frac{1}{2}\beta(\beta-1)%
    +\frac{1}{2}\delta(\delta-1),\\
B &= N-m,\\
C &= (m-mN+\frac{1}{2}N(N-1))+(m+1-N)\beta+(m-1)\delta.
\end{align*}
This can be written as a product of two sums which
simplifies to zero.
\end{proof}

Define a bilinear form
$(\ |\ ) : U^{\geq0}\times U^{\leq0} \longrightarrow \C(q)$ by
\begin{equation}\label{defbf}
(x|y) = \zeta(y)(x) \qquad(x\in U^{\geq0}, y\in U^{\leq0}).
\end{equation}

For $n\in\Z_{>0}$, we denote by
$\Delta_n : U \longrightarrow U^{\otimes(n+1)}$,
the algebra homomorphism defined by
$\Delta_1=\Delta$, $\Delta_n = (\Delta\otimes1)\circ\Delta_{n-1}$,
and we write
\begin{equation}
\Delta_n(x) = \sum_{(x)_n} x_{(0)}\otimes x_{(1)}\otimes\cdots\otimes x_{(n)}.
\end{equation}

For homogeneous elements $x_i \in U^{\geq0}_{\beta_i}$,
$y_i \in U^{\leq0}_{-\gamma_i}$ ($i=1,2$), we define
$(x_1\otimes x_2|y_1\otimes y_2) = %
\theta(\beta_2,-\gamma_1)(x_1|y_1)(x_2|y_2)$
and extend it by linearity.
For $x\in U_\beta$, $y\in U_\gamma$, we will write
$\theta(x,y)$ to mean $\theta(\beta,\gamma)$ and define
$P : U\otimes U \longrightarrow U\otimes U$ by 
$P(x\otimes y) = \theta(x,y)y\otimes x$ on homogeneous elements
and extend it by linearity.

\begin{prop}\label{prfrm}
The bilinear form $(\ |\ )$ on $U^{\geq0}\times U^{\leq0}$ defined by
\eqref{defbf} satisfies:
\begin{alignat}{2}
(x|y_1y_2) &= (\Delta(x)|y_1\otimes y_2)&
			&(x\in U^{\geq0}, y_1,y_2\in U^{\leq0}),\\
(x_1x_2|y) &= (P(x_1\otimes x_2)|\Delta(y))&\quad
			&(x_1,x_2\in U^{\geq0}, y\in U^{\leq0}),\label{ntrv}\\
(q^h|q^{h'}) &= q^{-(h|h')}& &(h,h'\in P^\vee),\\
(q^h|\fik) &= 0,\\
(\eik|q^h) &= 0,\\
(\eik|\fjl) &= -\frac{1}{\xi_i}\delta_{i,j}\delta_{k,l}
\end{alignat}
for $i,j\in\I$, $k=1,2,\cdots,m_j$.\\
Moreover, the bilinear form on $U^{\geq0}\times U^{\leq0}$ satisfying
the above equations is unique.
\end{prop}
\begin{proof}
Everything including uniqueness is straightforward except for \eqref{ntrv}.
It is proved by induction. Here we just show the induction part.
We suppress the summation signs for simplicity.
Assume
$(x^1x^2|y^i) = \theta(x^1,x^2)(x^2\otimes x^1|\Delta(y^i))$ for $i=1,2$.
Then,
{\allowdisplaybreaks\begin{align*}
(x^1x^2|y^1y^2)
&= \zeta(y^1)\zeta(y^2)(x^1x^2)\\
&= \zeta(y^1)\otimes\zeta(y^2)
((x^1_{(0)}\otimes x^1_{(1)})\cdot(x^2_{(0)}\otimes x^2_{(1)}))\\
&= \zeta(y^1)\otimes\zeta(y^2)
(\theta(x^1_{(1)},x^2_{(0)})x^1_{(0)}x^2_{(0)}\otimes x^1_{(1)}x^2_{(1)})\\
&= \theta(x^1_{(1)},x^2_{(0)})\theta(y^2,x^1_{(0)})\theta(y^2,x^2_{(0)})
(x^1_{(0)}x^2_{(0)}|y^1)(x^1_{(1)}x^2_{(1)}|y^2)\\
&= \theta(x^1_{(0)},x^1_{(1)})\theta(x^1_{(0)},x^2_{(1)})
\theta(x^2_{(0)},x^2_{(1)})\theta(x^1_{(0)},x^2_{(0)})
\theta(x^1_{(1)},x^2_{(1)})\\
&\phantom{=\theta(x}\cdot
(x^2_{(0)}\otimes x^1_{(0)}|y^1_{(0)}\otimes y^1_{(1)})
(x^2_{(1)}\otimes x^1_{(1)}|y^2_{(0)}\otimes y^2_{(1)})\\
&= \theta(x^1_{(0)},x^1_{(1)})\theta(x^1_{(0)},x^2_{(1)})
\theta(x^2_{(0)},x^2_{(1)})\\
&\phantom{=\theta(x}\cdot
(x^2_{(0)}|y^1_{(0)})(x^1_{(0)}|y^1_{(1)})
(x^2_{(1)}|y^2_{(0)})(x^1_{(1)}|y^2_{(1)})\\
&= \theta(y^1_{(1)},y^2_{(0)})
(\Delta(x^2)|y^1_{(0)}\otimes y^2_{(0)})
(\Delta(x^1)|y^1_{(1)}\otimes y^2_{(1)})\\
&= \theta(y^1_{(1)},y^2_{(0)})\theta(y^1_{(0)}y^2_{(0)},x^1)
(x^2\otimes x^1|y^1_{(0)}y^2_{(0)}\otimes y^1_{(1)}y^2_{(1)})\\
&= \theta(x^1,x^2)
(x^2\otimes x^1|%
(y^1_{(0)}\otimes y^1_{(1)})\cdot(y^2_{(0)}\otimes y^2_{(1)}))\\
&= \theta(x^1,x^2)(x^2\otimes x^1|\Delta(y^1y^2)).
\end{align*}}
This completes the proof.
\end{proof}

\begin{lem}\hfill
\begin{enumerate}
\item $(S(x)|S(y)) = (x|y)$
\item $(xq^h|yq^{h'}) = %
	q^{-(h|h')}(x|y) \quad(h,h'\in P^\vee, x\in U^+, y\in U^-)$.
\item $(U^+_\beta|U^-_{-\gamma}) = 0$ \quad if $\gamma\neq\beta$.
\end{enumerate}
\end{lem}
\begin{proof}
To prove (a), we set $(\ |\ )' = (S(\ )|S(\ ))$ and show $(\ |\ )'$
satisfies conditions of Proposition~\ref{prfrm}. The remaining two are easy.
\end{proof}

\begin{lem}\label{flip}
For $x\in U^{\geq0}$, $y\in U^{\leq0}$, that are homogeneous, we have
\begin{equation}\label{flip1}
\theta(x,y)yx =
\sum_{(x)_2,(y)_2}
\Theta_{xy}
(x_{(0)}|S(y_{(0)}))(x_{(2)}|y_{(2)})x_{(1)}y_{(1)}
\end{equation}
and
\begin{equation}\label{flip2}
xy = \hspace{-3mm}  \sum_{(x)_2,(y)_2} \hspace{-2mm}
                   \theta(x_{(1)},y_{(1)}) \Theta_{xy}
                   (x_{(0)}|y_{(0)})(x_{(2)}|S(y_{(2)})) y_{(1)}x_{(1)}
\end{equation}
with $\Theta_{xy} =
\theta(x_{(1)},y_{(0)})\theta(x_{(2)},y_{(0)})\theta(x_{(2)},y_{(1)}).$
\end{lem}
\begin{proof}
By substituting \eqref{flip1} into the right hand side of \eqref{flip2},
we can show that \eqref{flip1} implies \eqref{flip2}.

To prove \eqref{flip1}, we use induction on $y$ and reduce the problem to
showing this true for $y=q^h$ and $y=\fik$.
The case $y=q^h$ is easy. The case $y=\fik$ turns out to be equivalent to
showing
\begin{align}
\theta(x,\fik) \fik x =
\sum_{x_{(1)}}\Big\{
(x_{(1)}&|\fik)x_{(0)}\notag\\
&+ \theta(x_{(1)},\fik) (x_{(1)}|K_i^{-1})x_{(0)}\fik\\
&- \theta(x_{(1)},\fik) (x_{(0)}|\fik)K_i^{-1}x_{(1)} \notag \Big\},
\end{align}
which is proved by induction on the length of $x$.
\end{proof}

\begin{lem}\label{yzero}
Let $\beta\in Q^+\setminus\{0\}$ and $y\in U^-_{-\beta}$.
If $\eik y = \theta(\ali,-\beta)y\eik$
for all $i\in\I$, $k=1,2,\cdots,m_i$, then $y=0$.
\end{lem}
\begin{proof}
Choose $\lambda\in P^+$ satisfying the assumptions of
Lemma~\ref{UViso}.
Since $\eik(y\cdot v_\lambda) =
\theta(\ali,-\beta)y(\eik\cdot v_\lambda) = 0$
for all $i\in\I$, $k=1,2,\cdots,m_i$,
and $\text{wt}(y\cdot v_\lambda) = \lambda-\beta\lneq \lambda$,
$y\cdot v_\lambda$ generates a proper submodule of $V^q(\lambda)$.
Hence $y\cdot v_\lambda=0$.
Lemma~\ref{UViso} now says $y=0$.
\end{proof}

\begin{thm}\label{nondg}
For $\beta\in Q^+$, the bilinear form
$(\ |\ ) : U^{\geq0}\times U^{\leq0} \longrightarrow \C(q)$
defined by \eqref{defbf} is nondegenerate
when restricted to $U^+_\beta\times U^-_{-\beta}$.
\end{thm}
\begin{proof}
Since $\dim U^+_\beta = \dim  U^-_{-\beta}$, nondegeneracy on one side
implies the nondegeneracy on the other side.
So we will just prove the statement
\begin{equation}\label{eqb}
\text{if}\quad y\in U^-_{-\beta},\quad\text{and}
\quad ( U^+_\beta | y) = 0 \quad\text{then}\quad y=0.
\end{equation}
We use induction on $\beta$.\\
The case $\beta = 0\text{ or }\ali$ is easy.\\
Assume \eqref{eqb} is true for all $\gamma<\beta$ with
$\beta\in Q^+\setminus(\{0\}\cup\{\ali\}_{i\in\I})$.
Recall the notation $K_\gamma = q^{\sum_{i\in\I}n_i s_i h_i}$
for $\gamma=\sum_{i\in\I}n_i\ali \in Q$.
By definition of $\Delta$, we see that
\begin{equation}
\Delta(y) = \sum_{0\leq\gamma\leq\beta} y_\gamma(1\otimes K_{-\gamma})
,\quad y_\gamma\in U^-_{-\gamma}\otimes U^-_{-(\beta-\gamma)}
\end{equation}
with $y_0 = 1\otimes y$ and $y_\beta = y\otimes1$.
Fix $0<\gamma<\beta$. For any $u\in U^+_{\beta-\gamma}$ and
$v\in U^+_\gamma$, we have
\begin{align}
(v\otimes u|y_\gamma) &=
(v\otimes u|y_\gamma(1\otimes K_{-\gamma}))\\
&= (v\otimes u|\Delta(y))\\
&= \theta(\gamma,\beta-\gamma)(uv|y)\qquad\text{by~\eqref{ntrv}}\\
&= 0.
\end{align}
Hence $(U^+_\gamma\otimes U^+_{\beta-\gamma}|y_\gamma) = 0$.
This implies $y_\gamma = 0$ by our induction hypothesis.
Therefore $\Delta(y) = y\otimes K_{-\beta} + 1\otimes y$.
We apply Lemma~\ref{flip} to
\begin{align}
\Delta_2(\eik) &= \eik\otimes1\otimes1 + K_i\otimes\eik\otimes1
+ K_i\otimes K_i\otimes\eik\\
\Delta_2(y) &= y\otimes K_{-\beta}\otimes K_{-\beta} +
1\otimes y\otimes K_{-\beta} + 1\otimes1\otimes y
\end{align}
and get
\begin{equation}
\theta(\ali,-\beta)y\eik = \eik y \quad\text{for all}\quad i\in\I.
\end{equation}
Hence $y=0$ by Lemma~\ref{yzero}.
\end{proof}

\subsection{The Killing form}
Recall from Proposition~\ref{tridcmp} that
$U \cong U^+ \otimes U^0 \otimes S(U^-) \cong U^- \otimes U^0 \otimes S(U^+)$.
Using the bilinear form defined in the previous section, we define a new
bilinear form
\begin{equation*}
\langle\ | \ \rangle : U \times U \longrightarrow \C(q^{\frac{1}{2}})
\end{equation*}
by setting,
\begin{equation}
\langle x_1 q^{h_1} S(y_1) | y_2 q^{h_2} S(x_2) \rangle
= (x_1|y_2)(x_2|y_1)q^{-(h_1|h_2)/2} \theta(y_1,y_2)\theta(y_1,x_2)
\end{equation}
for homogeneous $x_i \in U^+$, $y_i \in U^-$, $h_i \in P^\vee$
and extending by linearity.

For homogeneous $u, v \in U$, we define
\begin{align}
\ad(u)\cdot v =& \sum_{(u)_1} \theta(u_{(1)}, v) u_{(0)} v S(u_{(1)})\\
v\cdot\adt(u) =& \sum_{(u)_1} \theta(v, u_{(0)})S(u_{0}) v u_{(1)}.
\end{align}
It is easy to check that these define left and right actions of $U$ on
$U$.

The bilinear form on $U$ defined above is \defi{invariant} in that:
\begin{prop}\label{action}
For $u, v, v' \in U$,
\begin{equation}
\langle \ad(u)\cdot v | v' \rangle = \langle v | v'\cdot\adt(u) \rangle
\theta(u,v) \theta(u, v').
\end{equation}
\end{prop}
\begin{proof}
It suffices to check the formula for $u= q^{h''}$ $(h'' \in P^\vee)$,
$\eik$, $\fik$ $(i \in\I, k=1,2, \cdots, m_i)$
and for $v=x q^{h} S(y)$
and $v'=y' q^{h'} S(x')$ with $x \in U^+_{\beta}$, $x' \in U^+_{\beta'}$,
$y \in U^-_{-\gamma}$, $y'\in U^-_{-\gamma'}$ $(\beta, \beta',
\gamma, \gamma'\in Q^+)$.
Since the case $u=\fik$ is similar to the case $u=\eik$, we will
omit the case $u=\fik$.

(i) $u= q^{h''}$\\
The left hand side is
\begin{equation*}
\langle \ad(u) \cdot v | v' \rangle
= \langle q^{h''} v q^{-h''} | v' \rangle
= q^{(\beta - \gamma) (h'')} \langle v | v' \rangle,
\end{equation*}
and the right hand side is
\begin{equation*}
\langle v | v' \cdot \adt(u) \rangle
= \langle v | q^{-h''} v' q^{h''} \rangle
= q^{(\gamma' - \beta')(h'')} \langle v | v' \rangle.
\end{equation*}
Since $\langle v | v' \rangle \neq 0$ only when $\beta = \gamma'$
and $\beta' = \gamma$, we are done.

(ii) $u= \eik$\\
Applying Lemma~\ref{flip}, we obtain
\begin{align*}
\ad (u) \cdot v & = \eik x q^{h} S(y) + \theta(\alpha_i, \beta)
q^{(\alpha_i | \beta)} x K_i q^h S(\eik y) \\
& = \eik x q^{h} S(y) + \theta(\alpha_i, \beta) q^{(\alpha_i | \beta)}
\sum _{(y)_2} \Big\{ A-B+C \Big\},
\end{align*}
where
\begin{align*}
A &= (\eik | y_{(0)})(1|S(y_{(2)})) x K_i q^{h} S(y_{(1)}),\\
B &= q^{\ali(h)} (K_i|y_{(0)}) (1|S(y_{(2)})) x\eik q^hS(y_{(1)}),\\
C &= \theta(\eik, y_{(1)}) (K_i|y_{(0)}) (\eik|S(y_{(2)}))x q^h S(y_{(1)}).
\end{align*}
and
\begin{align*}
v' \cdot \adt(u) &=-\theta(\beta'-\gamma',\alpha_i)q^{(\gamma'-\alpha_i|
\alpha_i)} \eik y'K_i^{-1} q^{h'}S(x') \\
& \qquad -\theta(\beta',\alpha_i)
q^{(\gamma'-\beta'|\alpha_i)}y'q^{h'}S(\eik x') \\
&=-\theta(\beta'-\gamma',\alpha_i)q^{(\gamma'-\alpha_i|\alpha_i)}
\sum_{(y')_2} \Big\{A'-B'+C' \Big\} \\
&\qquad -\theta(\beta', \alpha_i) q^{(\gamma'-\beta'|\alpha_i)}
y'q^{h'}S(\eik x'),
\end{align*}
where
\begin{align*}
A' &= (\eik|y'_{(0)})(1|S(y'_{(2)})) y'_{(1)} K_i^{-1}q^{h'}S(x'),\\
B' &= \theta(\ali,\beta'-\gamma') q^{(\alpha_i|\alpha_i)} q^{-\alpha_i(h')}
      (K_i|y'_{(0)}) (1|S(y'_{(2)})) y'_{(1)} q^{h'} S(x'\eik),\\
C' &= \theta(\eik, y'_{(1)}) (K_i|y'_{(0)}) (\eik|S(y'_{(2)}))
      y'_{(1)} q^{h'} S(x').
\end{align*}

There are only two cases to consider.
\begin{itemize}
\item $\gamma' = \beta + \alpha_i$ and $\gamma = \beta'$,
\item $\gamma'=\beta$ and $\gamma=\beta'+\alpha_i$.
\end{itemize}
Since the latter case is similar to the former, we will only check
the first case.
Assume $\gamma'=\beta+\alpha_i$ and $\gamma =\beta'$.
Then, in order to have $B\neq0$, we must have $y_{(0)}, y_{(2)} \in U^0$
and $y_{(1)} \in U^-_{-\gamma'}$.
Similarly, $A'\neq0$ implies $y'_{(0)} \in U^{\leq0}_{-\alpha_i}$,
$y'_{(1)} \in U^{\leq0}_{-\beta}$, and $y'_{(2)} \in U^0$.
In this case, we get $y'_{(1)} = \tilde{y}'_1 K_i^{-1}$ for some
$\tilde{y}'_1 \in U^{-}_{-\beta}$.
Also $C' \neq 0$ implies
$y'_{(0)} \in U^0$, $y'_{(1)} \in U^{\leq0}_{-\beta}$,
and $y'_{(2)} \in U^{\leq0}_{-\alpha_i}$.
In this case, we have $y'_{(2)} = \tilde{y}'_2 K_{\gamma' -\alpha_i}^{-1}$
for some $\tilde{y}'_2 \in U^{-}_{-\alpha_i}$.
We need to prepare one more fact.
Using Proposition~\ref{prfrm}, we obtain the following formula.
\begin{align*}
(x_1x_2x_3|y) = \sum_{(y)_2} &\theta(x_1x_2, x_3)
\theta(x_1x_2,y_{(0)})\theta(x_1,x_2)\theta(x_1, y_{(1)})\\
&\times(x_3|y_{(0)})(x_2|y_{(1)})(x_1|y_{(2)})
\end{align*}
for any $x_i\in U^+$ ($i=1,2,3$) and $y\in U^-$.
 From this formula, we get
\begin{align*}
&(x'|y) =(x'K_i | y)=\sum_{(y)_2}(K_i|y_{(0)})(x'|y_{(1)})(1|y_{(2)}), \\
&(x\eik|y')=\sum_{(y')_2}\theta(x,\eik)\theta(x,y'_{(0)})(\eik|y'_{(0)})
(x|y'_{(1)})(1|y'_{(2)}),\\
&(\eik x|y')=(\eik x K_i|y')\\
& \phantom{(\eik x|y')}=\sum_{(y')_2}
\theta(\eik,x)\theta(\eik,y'_{(1)})(K_i|y'_{(0)})(x|y'_{(1)})(\eik|y'_{(2)}).
\end{align*}

Now, we obtain
\begin{align*}
\langle \ad(u) &\cdot v|v' \rangle =\langle \eik xq^h S(y)|y'q^{h'}
S(x') \rangle \\
&\phantom{======} -\theta(\alpha_i,\beta)q^{(\alpha_i|\beta)}q^{\alpha_i(h)}
\sum_{(y)_2}(K_i|y_{(0)})
\langle x\eik q^h S(y_{(1)})|y' q^{h'}S(x') \rangle \\
& = \theta(\gamma, \gamma' -\beta') q^{-(h|h')/2} \\
&\phantom{==} \times \Big \{ (\eik x|y')(x'|y)
 -\theta(\alpha_i, \beta) q^{(\alpha_i|\beta)} q^{\alpha_i(h)}
\sum_{(y)_2} (K_i|y_{(0)} ) (x \eik|y')(x'|y_{(1)}) \Big \} \\
&= \theta(\gamma, \gamma'-\beta') q^{-(h|h')/2}\\
&\phantom{==} \times \sum_{(y)_2, (y')_2}
              (K_i|y_{(0)} ) (x'|y_{(1)} ) (x|y'_{(1)} )\\
&\phantom{======}\times\Big \{ (K_i|y'_{(0)} )( \eik|y'_{(2)} )
-\theta(x, y'_{(0)})q^{ (\alpha_i|\beta)} q^{\alpha_i(h)}
(\eik|y'_{(0)}) \Big \}
\end{align*}
and
{\allowdisplaybreaks
\begin{align*}
\langle v|v' \cdot \adt(u) \rangle =& -\theta(\beta'-\gamma', \alpha_i)
q^{(\gamma'-\alpha_i|\alpha_i)}\\
&\phantom{=}\times \sum_{(y')_2} \Big \{ (\eik|y'_{(0)})
\langle xq^h S(y)|y'_{(1)} K_i^{-1}q^{h'}S(x') \rangle \\
&\phantom{====} +\theta(\eik, y'_{(1)}) (K_i|y'_{(0)}) ( \eik|S(y'_{(2)}))
\langle xq^h S(y)|y'_{(1)} q^{h'}S(x') \rangle \Big \} \\
=& \theta(\beta' -\gamma', \alpha_i) \theta(\gamma,\beta-\beta')
q^{-(h|h')/2} q^{(\beta|\alpha_i)} \\
&\phantom{=}\times\sum_{(y')_2}
 \Big\{ \theta(\eik, y'_{(1)}) (K_i|y'_{(0)})(\eik K_i^{-1}|
\tilde{y}'_2 K_{\gamma'-\alpha_i}^{-1})(x|y'_{(1)})(x'|y)\\
&\phantom{====} - (\eik|y'_{(0)})( x|\tilde{y}'_1)(x'|y)
                  q^{\alpha_i(h)} \Big\} \\
=& \theta(\beta'-\gamma', \alpha_i) \theta(\gamma,\beta -\beta')
q^{-(h|h')/2}\\
&\phantom{=}\times \sum_{(y)_2,(y')_2} (K_i|y_{(0)})(x'|y_{(1)})(x|y'_{(1)})\\
&\phantom{====}\times\Big \{ \theta(\eik, y'_{(1)})
(K_i|y'_{(0)}) (\eik|y'_{(2)}) - q^{ (\alpha_i|\beta)} q^{\alpha_i(h)}
(\eik|y'_{(0)}) \Big \}.
\end{align*}
} 
Comparing these two, we get the desired formula.
\end{proof}

This proposition allows us to define a right $U$-module structure on some
subalgebra of $U^*$.
Define $\zeta : U \longrightarrow U^*$ by setting
\begin{equation}\label{zeta}
[\zeta(u)](v) = \langle v | u \rangle
\end{equation}
for $u, v \in U$.
Here, the dual space on the right should be viewed as the set of
linear maps from $U$ to $\C(q^\frac{1}{2})$.
For $\zeta(u) \in \zeta(U)$, $x\in U$, define $\zeta(u)\cdot x$ by,
\begin{equation*}
[\zeta(u)\cdot x](v) = \theta(u,x)\theta(v,x) [\zeta(u)](\ad(x)\cdot v).
\end{equation*}
Proposition~\ref{action} allows us to check
$\zeta(u)\cdot x = \zeta(u\cdot\adt(x))$.
So this gives a right $U$-module structure on $\zeta(U)$ and
$\zeta : U \longrightarrow \zeta(U)$ becomes a $U$-module homomorphism.

\begin{prop}\label{pr:48}
The bilinear form $\langle\ | \ \rangle$ is nondegenerate.
Hence, the map $\zeta$ is injective.
\end{prop}
\begin{proof}
Let $u\in U^-_{-\alpha}U^0 S(U^+_\beta)$ with $\langle v|u\rangle = 0$
for all $v\in U^+_\alpha U^0 S(U^-_{-\beta})$.
It suffices to show $u=0$.
For each $\gamma\in Q^+ - \{0\}$, choose a basis $\{u^\gamma_i\}_i$
of $U^+_\gamma$.
And let $\{v^\gamma_i\}_i$ be a basis of $U^-_{-\gamma}$ dual to 
$\{u^\gamma_i\}_i$ with respect
to the nondegenerate bilinear form $(\ |\ )$.
Notice that the elements $u^\alpha_i q^h S(v^\beta_j)$ with $h\in P^\vee$
and $i$, $j$ going over appropriate indices, form a basis for
$U^+_\alpha U^0 U^-_{-\beta}$.
Similarly, the elements $v^\alpha_k q^{h'} S(u^\beta_l)$ form a basis
for $U^-_{-\alpha} U^0 U^+_\beta$.
Writing $u = \sum_{k,h,l} a_{k,h,l} v^\alpha_k q^{h'} S(u^\beta_l)$
with $a_{k,h,l}\in\C(q)$,
and using
\begin{equation*}
\langle u^\alpha_i q^h S(v^\beta_j) |
        v^\alpha_k q^{h'} S(u^\beta_l) \rangle
= \delta_{i,k} \delta_{j,l} q^{-(h|h')/2}
  \theta(\beta,\alpha)\theta(\beta,\beta),
\end{equation*}
we arrive at,
\begin{equation*}
\sum_{h'\in P^\vee} a_{k,h',l} q^{-(h|h')/2} = 0
\end{equation*}
for each $k$, $l$, and $h\in P^\vee$.
Now, each map $h \mapsto q^{-(h|h')/2}$ is a group homomorphism from
$P^\vee$ to the multiplicative group $\C(q^\frac{1}{2})^\times$.
Since $q^{\frac{1}{2}}$ is not a root of unity, distinct $h'$ produces
distinct homomorphisms.
So, by Artin's Theorem on linear independence of characters,
every $a_{k,h',l} = 0$.
We have $u = 0$ as claimed.
\end{proof}

\section{Harish-Chandra Homomorphism}\label{HC}

We denote the center of $U$ by $\cntr$.
For each $i\in\I$ with $a_{i,i}\neq 0$, define the \defi{simple reflection}
$r_i \in {\rm GL}(\h)$ by,
\begin{equation*}
r_i(h) = h - \frac{2}{a_{i,i}}\ali(h)h_i,
\end{equation*}
and let
$\widetilde{W} = \langle r_i\; |\; i\in \I, a_{i,i} \neq 0 \rangle%
\subset {\rm GL}(\h)$.
Let $(U^0)^{\widetilde{W}}$ be the subspace of $U^0$ consisting of the
elements $\sum_{h\in P^\vee} c_h q^h$ ($c_h\in\C(q)$) such that $c_h\neq0$
implies $w(h)\in P^\vee$ and $c_{w(h)} = c_h$ for any $w\in \widetilde{W}$.

We define an algebra automorphism $\phi:U^0\rightarrow U^0$ by setting
$\phi(q^h) = q^{-\rho(h)}q^h$ for $h\in P^\vee$.
The \defi{Harish-Chandra homomorphism} $\xi : \cntr\rightarrow U^0$
is the restriction to $\cntr$ of the map
\begin{equation*}
U \xrightarrow{_\sim} U^-\otimes U^0 \otimes U^+
\xrightarrow{\varepsilon\otimes1\otimes\varepsilon} U^0
\xrightarrow{\phi} U^0.
\end{equation*}
For later use, we define the algebra homomorphism
$\chi_\lambda : U^0 \longrightarrow \C(q)$ for each $\lambda\in P^+$ by
$\chi_\lambda(q^h) = q^{\lambda(h)}$.

\begin{prop}\hfill
\begin{enumerate}
\item $\xi$ is an algebra homomorphism.
\item $\xi$ is injective. \label{injective}
\end{enumerate}
\end{prop}
\begin{proof}
We will just prove \ref{injective}.
Let $z\in\cntr$ be such that $\xi(z)=0$. Writing
$z=\sum_{\beta\in Q^+}z_\beta$
with $z_\beta \in U^-_{-\beta} U^0 U^+_\beta$, we see that $z_0 = 0$.
Fix any $\beta\in Q^+$ minimal with the property that $z_\beta \neq 0$.
Also choose basis $\{y_r\}_r$ and $\{x_s\}_s$ of
$U^-_{-\beta}$ and $U^+_\beta$, respectively.
We may write $z_\beta = \sum_{r,s} y_r u_{r,s} x_s$ for some $u_{r,s}\in U^0$.
Then,
\begin{align*}
0 = \eik z - z \eik =
&\sum_{\gamma\neq\beta}(\eik z_\gamma - z_\gamma \eik)\\
&+ \sum_{r,s} (\eik y_r -\theta(\ali,-\beta)y_r\eik)u_{r,s}x_s\\
&+ \sum_{r,s} y_r (\theta(\ali,-\beta)\eik u_{r,s}x_s - u_{r,s}x_s\eik)
\end{align*}
Recalling the minimality of $\beta$, we see that only the second term
on the right belongs to $U^-_{-(\gamma-\ali)}U^0 U^+_\gamma$.
So we have
$\sum_{r,s}(\eik y_r - \theta(\ali,-\beta)y_r\eik)u_{r,s}x_s = 0$.
$\{x_s\}_s$ was chosen to be a basis, so
$\eik\sum_r y_r u_{r,s} = \theta(\ali,-\beta)\sum_r y_r\eik u_{r,s}$
for all $i\in\I$ and $s$.

Let $v_\lambda\in V^q(\lambda)$
denote the highest weight vector.
Set $v = \sum_r \chi_\lambda(u_{r,s})y_r v_\lambda$.
Then $\eik v = \theta(\ali,-\beta)\sum_r y_r\eik u_{r,s} v_\lambda = 0$
for all $i\in\I$, so the irreducibility of $V^q(\lambda)$ says $v=0$.
Choosing an appropriate $\lambda\in P^+$, we may use Lemma~\ref{UViso}
and say $\sum_r \chi_\lambda(u_{r,s})y_r = 0$. Again, $\{y_r\}_r$ was a basis,
so $\chi_\lambda(u_{r,s}) = 0$ for all $r,s$. By choosing a suitable
set of $\lambda$, we may show $u_{r,s}=0$ for all
$r,s$ and we have $z_\beta=0$. This contradicts the choice of $z_\beta$.
\end{proof}

We now try to close in on the image of $\xi$.
For each $J\subset \{ (i,k)\;|\; i\in\I, k=1,2,\cdots,m_i\}$, let
$U_J = \langle \eik, \fik, U^0 \;|\; (i,k)\in J \rangle$.
We denote by $\cntr_J$ the center of the algebra $U_J$ and
by $\xi_J : \cntr_J \longrightarrow U^0$ the Harish-Chandra homomorphism
for $U_J$.
Let $U^+_J$ (respectively, $U^-_J$) be the subalgebra of $U_J$
generated by $\eik$ (respectively $\fik$) with $(i,k) \in J$, and set
\begin{alignat}{2}
R^+_J &= \{ x\in U^+ |\; (x|U^-_J) = 0 \}
      &= \{ x\in U^+ |\; (x|U^-_J U^0) = 0 \}\\
R^-_J &= \{ y\in U^- |\; (U^+_J|y) = 0 \}
      &= \{ y\in U^- |\; (U^0 U^+_J|y) = 0 \}\\
R_J &= R^-_J U^0 U^+ + U^- U^0 R^+_J
\end{alignat}
The following may be proved as in \cite{KgTa}.
\begin{lem}\label{sum}\hfill
\begin{enumerate}
\item $U = U_J \oplus R_J$,
\item $U_J R_J U_J \subset R_J$,
\item $(\varepsilon\otimes 1 \otimes\varepsilon)(R_J) = 0$.
\end{enumerate}
\end{lem}

Define $U^0_r = \bigoplus_h \C(q) q^h$, where the direct sum is over all
$h \in P^\vee$ satisfying,
\begin{itemize}
\item $\ali(h) \in s_i a_{i,i} \Z$ if $i \in \I^{\text{ev}}$,
\item $\ali(h) \in 2s_i a_{i,i} \Z$ if $i \in \I^{\text{odd}}$
and $a_{i,i}\neq0$.
\end{itemize}

\begin{prop}\label{prop:53}\hfill
\begin{enumerate}
\item $\im(\xi) \subset (U^0)^{\widetilde{W}}$. \label{Weyl-inv}
\item $\im(\xi) \subset U^0_r$. \label{si}
\item $\im(\xi) \subset \im(\xi_J)$. \label{xiJ}
\end{enumerate}
\end{prop}
\begin{proof}
\ref{Weyl-inv} Let $z \in \cntr$.
Let $v_\lambda\in M^q(\lambda)$ be
the highest weight vector. Then,
$z v_\lambda = \chi_{\lambda+\rho}(\xi(z))v_\lambda$.
Since $z$ commutes with every element of $U$, $z$ acts as
$\chi_{\lambda+\rho}(\xi(z))$ on every element of $M^q(\lambda)$.
Now, fix $i\in\I$ such that $a_{i,i}\neq0$. We may calculate
\begin{equation*}
\eik\fik^n = \theta_{i,i}^n\fik^n\eik +
\theta_{i,i}^{n-1}\fik^{n-1}\frac{1}{\xi_i}
\Big(\frac{1-\theta_{i,i}^{-n}q_i^{-na_{i,i}}}
          {1-\theta_{i,i}^{-1}q_i^{-a_{i,i}}} K_i -
     \frac{1-\theta_{i,i}^nq_i^{na_{i,i}}}
          {1-\theta_{i,i}q_i^{a_{i,i}}} K_i^{-1}\Big).
\end{equation*}
So that, for each $\lambda\in P$ satisfying
$n(\lambda) := \frac{2}{a_{i,i}}\lambda(h_i) \in \Z_{\geq0}$,
we can check that $\fik^{n(\lambda)+1}v_\lambda$
is a highest weight vector.
Its weight is
\begin{align*}
\lambda - \Big(\frac{2}{a_{i,i}}\lambda(h_i) + 1 \Big)\ali
& = \lambda - \frac{2}{a_{i,i}}(\lambda+\rho)(h_i)\ali\\
& = r_i(\lambda+\rho) - \rho.
\end{align*}
The argument at the beginning of this proof applies to any highest weight
vector and we have,
\begin{equation*}
\chi_{\lambda+\rho}(\xi(z)) = \chi_{r_i(\lambda+\rho)}(\xi(z))
\end{equation*}
under the condition $\frac{2}{a_{i,i}}\lambda(h_i) \in \Z$.
Checking $\chi_{r_i\mu}(q^h) = \chi_\mu(r_i q^h)$, for any $\mu\in\h^*$,
the above may now be written as
\begin{equation*}
\chi_{\lambda+\rho}(\xi(z) - r_i \xi(z)) = 0
\end{equation*}
for every $\lambda\in P$ satisfying
$\lambda(h_i)\in \frac{a_{i,i}}{2}\Z_{\geq0}$.
By choosing a suitable set of $\lambda$, we may show $\xi(z) = r_i\xi(z)$.

\ref{si} Let $z= \sum_{\beta\in Q^+} z_\beta \in \cntr$ with
$z_\beta \in U^-_{-\beta} U^0 U^+_\beta$.
Set $x = \sum_{n=0}^{\infty} z_{n\ali}$ and
$y = z - x$.
Then $z = x + y$ with $x\in U_{\{(i,1)\}}$ and $y\in R_{\{(i,1)\}}$.
Looking at
\begin{equation*}
0 = \eik z - z\eik = (\eik x - x\eik) + (\eik y - y\eik)
\end{equation*}
with Lemma~\ref{sum} in mind, we see that $x\in\cntr_{\{(i,1)\}}$.
By the results of Section~\ref{rank1}, all of which may be obtained by
direct calculation, we have,
\begin{itemize}
\item $z_0\in\langle K_i,q^h\;|\;\ali(h)=0\rangle$ if $i\in\I^{\text{ev}}$,
\item $z_0\in\langle K_i^2,q^h\;|\;\ali(h)=0\rangle$ if $i\in\I^{\text{odd}}$
      and $a_{i,i}\neq0$.
\end{itemize}
The result follows.

\ref{xiJ} For $z \in \cntr$, write $z = x + y$ with $x\in U_J$
and $y\in R_J$. As in the proof for \ref{si}, we may show $x\in\cntr_J$.
So we have,
$\xi(z) = \xi(x) + \xi(y) = \xi(x) = \xi_J(x) \in \im(\xi_J)$.
\end{proof}

\section{The Center of $U_q(\g)$}\label{CU}

\subsection{Rank 1}\label{rank1}
In this section, we list the center for the case when the index set is of
size 1. All results may be obtained by direct calculation using induction
after choosing a suitable basis of $U\cong U^-\otimes U^0\otimes U^+$.\\
If $a_{i,i} \neq 0$, $\theta_{i,i} = 1$, define
\begin{equation*}
C_i =  f_{i,1}e_{i,1} + \frac{1}{\xi_i}
\Big(
\frac{1}{1-q^{-s_i a_{i,i}}}K_i - \frac{1}{1-q^{s_i a_{i,i}}}K_i^{-1}
\Big).
\end{equation*}
If $a_{i,i} \neq 0$, $\theta_{i,i} = -1$, define
\begin{align*}
C_i = f_{i,1}^2 & e_{i,1}^2
+
\frac{1}{\xi_i} f_{i,1}
\Big(
\frac{1-q^{s_i a_{i,i}}}{1+q^{s_i a_{i,i}}}K_i -
\frac{1-q^{-s_i a_{i,i}}}{1+q^{-s_i a_{i,i}}}K_i^{-1}
\Big)
e_{i,1}\\
&-
\frac{1}{\xi_i^2}\Big\{
\frac{1}{(1+q^{-s_i a_{i,i}})^2}K_i^2 +
\frac{1}{(1+q^{s_i a_{i,i}})^2}K_i^{-2}
\Big\}.
\end{align*}
If $h\in P^\vee$ satisfy $\ali(h)\neq 0$, define
\begin{equation*}
C_{ih} = 
f_{i,1}q^h e_{i,1}+\frac{1}{\xi_i}\frac{1}{1-q^{-\ali(h)}}q^h(K_i - K_i^{-1}).
\end{equation*}

%
%

\begin{prop}\label{6.2}\hfill
\begin{enumerate}
\item If $J = \{(i,1)\}$ and $a_{i,i} \neq 0$,\hfill\\ then
      $\cntr_J = \langle C_i, q^h \;|\; h\in P^\vee, \ali(h) = 0 \rangle$.
\item If $J= \{(i,1)\}$, $a_{i,i} =0$, and $\theta_{i,i} = 1$,\hfill\\ then
      $\cntr_J = \langle q^h \;|\;
      h\in P^\vee, \ali(h) = 0 \rangle \subset U^0$.\label{gau}
\item If $J= \{(i,1)\}$, $a_{i,i} =0$, and $\theta_{i,i} = -1$,\hfill\\ then
      $\cntr_J = \langle C_{ih}, q^{h'} \;|\;
      h, h' \in P^\vee, \ali(h) \neq 0, \ali(h') = 0 \rangle$.
\end{enumerate}
\end{prop}

\subsection{Finite type}
In this section, we give a structure theorem for the center of $U_q(\g)$
when the Borcherds-Cartan matrix is of finite type.
We take the Borcherds-Cartan matrix to be of finite type throughout
this section.
To simplify arguments, we redefine
\begin{align*}
P^\vee &= \bigoplus_{i\in\I}\Z h_i,\\
\h &= \bigoplus_{i\in\I} \C h_i
\end{align*}
for this section.
Notice that the bilinear form $(\ |\ )$ is still nondegenerate on
the redefined $\h$.

The irreducible highest weight module has a natural grading.
\begin{equation*}
V^q(\lambda) = \bigoplus_{\alpha\in Q^+} V^q(\lambda)_{\lambda-\alpha}.
\end{equation*}
Define a map $\eta \in \operatorname{End}(V^q(\lambda))$ by setting
$\eta(v) = \theta(\alpha,\alpha) v$ for $v\in V^q(\lambda)_{\lambda-\alpha}$.
When the Borcherds-Cartan matrix $A$ is of finite type,
it is known(\cite{Kac77}) that the irreducible highest weight module
$V(\lambda)$ over $\g(A)$ is finite dimensional for $\lambda\in P^+$.
Since the classical limit(\cite{Kang96}) of $V^q(\lambda)$ is
$V(\lambda)$, $V^q(\lambda)$ is also of finite dimension when
$\lambda\in P^+$.
So we may define the \defi{supertrace} for $x\in U_q(\g)$ acting on
$V^q(\lambda)$ by
\begin{equation}
\str(x ; V^q(\lambda)) = \tr(\eta\circ x ; V^q(\lambda)).
\end{equation}
For homogeneous elements $x,y \in U$, we can easily check
\begin{equation}\label{strflip}
\str(xy)=\theta(x,y)\str(yx).
\end{equation}

\begin{lem}\label{lm:64}
$u\in\cntr$ if and only if $u\cdot\adt(x) = \varepsilon(x) u$ for all
$x \in U$.
\end{lem}
\begin{proof}
Let $u\in\cntr$.
Then, $u\in U_0$ and
\begin{align*}
u\cdot\adt(x)
&= \sum_{(x)_1} S(x_{(0)}) u x_{(1)}\\
&= u \sum_{(x)_1} S(x_{(0)}) x_{(1)} = \varepsilon(x) u.
\end{align*}
Conversely, if $u\cdot\adt(x) = \varepsilon(x) u$ for all $x \in U$,
\begin{equation*}
q^{-h} u q^h = u \cdot \adt(x) = \varepsilon(q^h)u = u.
\end{equation*}
So $u\in U_0$ and we have,
\begin{equation*}
0 = \varepsilon(\eik)u = u\cdot\adt(\eik)
  = -K_i^{-1}\eik u + K_i^{-1} u\eik.
\end{equation*}
This shows $\eik u = u \eik$.
We may similarly show $\fik u = u \fik$ and hence $u\in\cntr$.
\end{proof}

For each $\lambda\in P^+$, define $f_{\lambda} \in U^*$ by
\begin{equation}
f_{\lambda}(u) = \str(uK_{2\rho}^{-1};V^q(\lambda)).
\end{equation}
Let $\nu : \h \rightarrow \h^*$ denote the isomorphism given by the
nondegenerate symmetric bilinear form $(\ |\ )$.
Define
\begin{equation*}
\widehat{Q} = \nu(P^\vee) = \bigoplus_{i\in\I} \Z\frac{1}{s_i}\ali.
\end{equation*}
Recall the map $\zeta: U \longrightarrow U^*$ defined in \eqref{zeta}.

\begin{lem}\label{lem:65}
For $\lambda\in P^+$, $f_\lambda\in \im(\zeta)$ if and only if
$\lambda\in\frac{1}{2}\widehat{Q}$.
\end{lem}
\begin{proof}
 From Proposition~\ref{pr:48}, we see that the image of $\zeta$ is the
restricted dual of $U_q(\g)$.
So
\begin{equation*}
\im(\zeta) =
(\bigoplus_{\beta\in Q^+} (U^-_{-\beta})^*)\otimes
(\bigoplus_{\mu\in\frac{1}{2}\widehat{Q}} \C(q)\chi_\mu)\otimes
(\bigoplus_{\beta\in Q^+} (U^+_\beta)^*)
\end{equation*}
under the identification $U\cong U^-\otimes U^0 \otimes U^+$.
The finite dimensionality of $V^q(\lambda)$ allows us to shows
$f_\lambda\in\im(\zeta)$ if and only if
$\lambda\in\frac{1}{2}\widehat{Q}$.
\end{proof}

The next proposition gives elements of the center.

\begin{prop}\label{pr:66}
For each $\lambda\in P^+\cap\frac{1}{2}\widehat{Q}$,
we have $z_\lambda := \zeta^{-1}(f_{\lambda}) \in \cntr$.
\end{prop}
\begin{proof}
Recall from the theory of finite dimensional simple Lie algebras, that
$\rho$ may be written as a half sum of positive roots.
Since the simple roots for the super case is identical to the non-super case,
we have $2\rho \in Q^+$ in either case.
Hence, in the notation given on page~\pageref{nata},
$K_{2\rho}$ is a well-defined element of $U^0$.
Using the fact that $K_{2\rho}^{-1}xK_{2\rho} = S^2(x)$ for any $x \in U$
and using the property of supertrace given by \eqref{strflip}, we have for
any $u\in U$,
{\allowdisplaybreaks
\begin{align*}
(f_{\lambda}\cdot x)(u) & = f_{\lambda}(\ad(x)\cdot u) \theta(u,x) \\
& = \sum_{(x)_1} \str(x_{(0)}uS(x_{(1)})K_{2\rho}^{-1};
V^q(\lambda)) \theta(x_{(1)}, u)\theta(u,x) \\
& = \sum_{(x)_1} \str(uS(x_{(1)})K_{2\rho}^{-1}x_{(0)};
V^q(\lambda) )\theta(x_{(0)}, x_{(1)}) \\
& = \str(uS(\sum_{(x)_1}S(x_{(0)})x_{(1)})K_{2\rho}^{-1};
V^q(\lambda)) \\
&= \varepsilon(x) \str(u K_{2\rho}^{-1};V^q(\lambda)) \\
&= \varepsilon(x) f_{\lambda}(u).
\end{align*}
}
Thus $f_{\lambda} \cdot x = \varepsilon(x) f_{\lambda}$.
Recall from Proposition~\ref{pr:48} that $\zeta$ is injective, and notice
\begin{equation*}
f_\lambda\cdot x = \zeta( \zeta^{-1}(f_\lambda))\cdot x
= \zeta( \zeta^{-1}(f_\lambda) \cdot \adt(x)).
\end{equation*}
This shows $\zeta^{-1}(f_{\lambda}) \cdot \adt(x)
= \varepsilon(x)\zeta^{-1}(f_{\lambda})$.
 From Lemma~\ref{lm:64}, we get $\zeta^{-1}(f_{\lambda}) \in \cntr$.
\end{proof}

We finally show that the above elements generate the whole center.

\begin{thm}
Suppose that the Borcherds-Cartan matrix $A = (a_{i,j})_{i,j\in\I}$
is indecomposable and of finite type.
Then, $\xi : \cntr \longrightarrow (U_r^0)^{\widetilde{W}}$ is an
isomorphism.
\end{thm}
\begin{proof}
Let us calculate $\xi(z_\lambda)$.
We extend the notation $K_\beta$ introduced in page~\pageref{nata}
to $\beta\in\widehat{Q}$ by setting $K_{\frac{1}{s_i}\ali} = q^{h_i}$.
We have the commutative diagram
\begin{center}
\begin{texdraw}
\textref h:C v:C
\arrowheadsize l:2.4 w:1.1 \arrowheadtype t:F
\htext(0 0){$U$}
\htext(21 0){$U^*$}
\htext(0 -13){$U^0$}
\htext(21 -13){$(U^0)^*$}
\move(4 0)\avec(16 0)
\move(4 -13)\avec(16 -13)
\move(0 -3)\avec(0 -10)
\move(20 -3)\avec(20 -10)
\htext(10 3){$\zeta$}
\htext(-8 -6){$\varepsilon\otimes\id\otimes\varepsilon$}
\end{texdraw}
\end{center}
where the right vertical arrow is the restriction map and the lower
horizontal arrow is given by $K_\mu \mapsto \chi_{-\mu/2}$.
Now, as maps on $U^0$,
\begin{equation*}
f_\lambda = \sum_{\mu\leq\lambda} \theta(\lambda-\mu,\lambda-\mu)
        \dim(V(\lambda)_\mu)q^{-2(\rho|\mu)}\chi_\mu.
\end{equation*}
This shows,
\begin{equation}\label{tsho}
\xi(z_\lambda) = \sum_{\mu\leq\lambda}
      \theta(\lambda-\mu,\lambda-\mu)\dim(V(\lambda)_\mu)K_{-2\mu}
\end{equation}
for $\lambda\in P^+\cap\frac{1}{2}\widehat{Q}$.

Define $\widehat{P}$ to be the set of elements $\mu\in\h^*$ such that
$\mu(h_i)\in\Z$ if $i\in\I$ is even and $\mu(h_i)\in2\Z$ if $i\in\I$
is odd.
Notice $P^+\subset\widehat{P}$.
We can now write
\begin{equation*}
U_r^0 = \bigoplus_{\mu\in 2\widehat{P}\cap\widehat{Q}} \C(q)K_\mu.
\end{equation*}
Action of the Weyl groups $\widetilde{W}$ and $W$ defined on $\h$ and $\h^*$
are compatible with the isomorphism $\nu$.
By Proposition~\ref{pr:66}, it suffices to show that the elements
$\xi(z_\lambda)$ with $\lambda\in P^+\cap\frac{1}{2}\widehat{Q}$
generate $(U_r^0)^W$.

Set $\bar{\mu} = \sum_{w\in W}K_{-w\mu}$ for any $\mu\in\widehat{Q}$.
We know that the elements $\bar{\mu}$ with $\mu\in 2P^+\cap\widehat{Q}$
generate $(U^0_r)^W$.
Let us use induction to show that each of them belong to $\im(\xi)$.
The element $\bar{0}\in U^0_r$ is given by $\xi(z_\lambda)$ with $\lambda=0$.
Choose any $\lambda\in 2P^+\cap\widehat{Q}$.
Then, $\frac{1}{2}\lambda\in P^+\cap\frac{1}{2}\widehat{Q}$
so that $z_{\frac{1}{2}\lambda}$ is an element of the center.
Recall that $\im(\xi)$ is invariant under the action of $W$
(Proposition~\ref{prop:53}).
Using $\dim V^q(\frac{1}{2}\lambda)_{\frac{1}{2}\lambda} = 1$, we may rewrite
\eqref{tsho} as
\begin{equation*}
\xi(z_{\frac{1}{2}\lambda}) = \bar{\lambda} + \sum_\mu n_\mu \bar{\mu},
\end{equation*}
with $n_\mu\in\Z$ and $\frac{1}{2}\mu$ running over some set of weights of
$V^q(\frac{1}{2}\lambda)$.
Since all $\mu\lneq\lambda$, induction hypothesis show that each $\bar{\mu}$
belong to $\im(\xi)$.
Hence $\bar{\lambda}\in\im(\xi)$ and the induction step is complete.
\end{proof}

\subsection{Other Cases}
Let $2_i$, $0_i$, and $\circleddash_i$ denote the fact that $a_{i,i}$ is 
respectively, 2, 0, and negative.
We will sometimes add a $\pm$ to these to reflect the sign of
$\theta_{i,i}$. So, for example, $2_i^-$ implies that $i$ is an odd real index.
For $i,j\in \I$, let us say $\odot_i$ is \defi{connected directly} to
$\odot_j$ if
$a_{i,j}\neq0$, where $\odot$ can be any one of 2, 0, or $\circleddash$.
Here are some results for the case when $|J|=2$.

\begin{lem}\label{rank2}
Assume one of the following.
\begin{enumerate}
\item $J=\{ (i,1),(i,2) \}$ with $\circleddash_i$ \label{ra}
\item $J=\{ (i,1),(j,1) \}$ with $0_i^+$ connected directly to \label{pm}
   $0_j^-$
\item $J=\{ (i,1),(j,1) \}$ with $0_i^-$ connected directly to \label{rf}
   $0_j^-$
\item $J=\{ (i,1),(j,1) \}$ with $\circleddash_i$ connected directly to
   $0_j^-$ \label{re}
\item $J=\{ (i,1),(j,1) \}$ with $2_i$ connected directly to $\circleddash_j$
   \label{rc}
\item $J=\{ (i,1),(j,1) \}$ with $\circleddash_i$ connected directly to
   $\circleddash_j$ \label{rb}
\end{enumerate}
Then, $\cntr_J \subset U^0$.
\end{lem}
\begin{proof}
\ref{ra} and \ref{rb} may be proved as in \cite[Proposition 4.5]{KgTa}.
And \ref{rc} may be proved as in \cite[Proposition 4.6]{KgTa}.
\ref{rf} is proved by explicit calculation.

Let us prove \ref{pm} and \ref{re} simultaneously.
Let $z\in \cntr_J$. Since it commutes with $q^h$ for all $h\in P^\vee$,
$z = \sum z_\beta$ with
$z_\beta \in U_\beta^- \otimes U^0 \otimes U^+_\beta$,
where the sum is over all
$\beta \in \Z_{\geq 0}\ali \oplus \Z_{\geq 0}\alj$.
Let $\alpha$ be maximal among those
$\beta \in \Z_{\geq 0}\ali \oplus \Z_{\geq 0}\alj$
for which $z_\beta$ is nonzero and suppose $\alpha\neq 0$.
Let $\{x_\mu\}$ and $\{y_\lambda\}$ be any bases of
$(U_J^+)_\alpha$ and $(U_J^-)_{-\alpha}$ respectively.
We can now write
\begin{equation*}
z = \Big(
\sum_{\lambda,\mu,h} c_h^{\lambda,\mu} y_\lambda q^h x_\mu
    \Big) + z'.
\end{equation*}
Recall Lemma~\ref{flip} and notice
\begin{align*}
\Delta_2(\eik) &= \eik\otimes1\otimes1 + K_i \otimes\eik\otimes1
               + K_i \otimes K_i \otimes\eik,\\
\Delta_2(y_\lambda) & = 1 \otimes y_\lambda \otimes K_{-\alpha}
               + \text{``other terms''}.
\end{align*}
This shows that the only part of $\eik z - z \eik$ belonging to the
direct sum component
 $U^-_{-\alpha}\otimes U^0 \otimes U^+_{\alpha+\ali}$ is
\begin{align*}
\Big(\sum_{\lambda,\mu,h} c_h^{\lambda,\mu} \theta(\ali,-\alpha)&
     y_\lambda \eik q^h x_\mu \Big) -
\Big(\sum_{\lambda,\mu,h} c_h^{\lambda,\mu} y_\lambda q^h x_\mu \eik \Big)\\
&= \sum_{\lambda,h} y_\lambda q^h \sum_{\mu}
c_h^{\lambda,\mu}
\Big( q^{-\ali(h)}\theta(\alpha,\ali) \eik x_\mu - x_\mu \eik \Big).
\end{align*}
Hence, for each $h\in P^\vee$ and $\lambda$,
\begin{equation*}
\sum_{\mu} c_h^{\lambda,\mu} \Big(
q^{-\ali(h)}\theta(\alpha,\ali) \eik x_\mu - x_\mu \eik \Big)
=0,
\end{equation*}
and the same statement with $i$ replaced by $j$ also holds.
Now, $e_{j,1}^2 = 0$ is the only relation in $U_J^+$ for the case we are
considering, so we may take an explicit set of monomials in $e_{i,1}$ and
$e_{j,1}$ for the basis of $U^+_\alpha$ and using these,
we can show that the two equations cannot be simultaneously true.
\end{proof}

\begin{prop}\label{mulguisin}
Assume that $A$ is indecomposable.
Suppose that every $0_j^-$ is connected directly to a
$0_i$ or a $\circleddash_i$.
If there is a nonempty subset $J$ of
$\{ (i,k) \;|\; i\in\I, k=1,\cdots,m_i \}$ such that $\cntr_J \subset U^0$,
then $\cntr$ is contained in $U^0$.
\end{prop}
\begin{proof}
Let $\bar{J} = \{ i\in\I\;|\; (i,k)\in J \text{ for some } k \}$.
For $i\in\I$, set
\begin{equation*}
T_i = \bigoplus_{h\in P^\vee, \ali(h)=0}\C(q)q^h.
\end{equation*}
We then have $\cntr \cap U^0 = \cap_{i\in\I} T_i$ and similarly,
$\cntr_J \cap U^0 = \cap_{i\in \bar{J}} T_i$.
It suffices to show $\im(\xi) \subset \cap_{i\in\I} T_i$.

We already have $\im(\xi) \subset \im(\xi_J) \subset T_i$
for every $i\in \bar{J}$.
Also if $0_i^+$, we have,
$\im(\xi)\subset\im(\xi_{\{(i,1)\}})\subset T_i$
by Proposition~\ref{6.2}~\ref{gau}.
If $0_i^-$, the conditions on the matrix shows we may use Lemma~\ref{rank2}
and say $\im(\xi)\subset T_i$.

We now show that if $a_{j,j}\neq0$ and $a_{i,j}\neq0$,
then $T_i \cap (U^0)^{\widetilde{W}} \subset T_j$.
Let $c = \sum c_h q^h \in T_i \cap (U^0)^{\widetilde{W}}$. We must have
$r_j c = c \in T_i$, so if $c_h\neq0$, then $\ali(h)=0$ and $\ali(r_j h)=0$.
But $\ali(r_j h) = - \frac{2}{a_{j,j}}a_{j,i}\alj(h)$ so $\alj(h)=0$.
We have $c\in T_j$ as wanted.

Fix any $j\in \I - (\bar{J}\cup\{i\;|\; a_{i,i}=0\})$. By the
indecomposability of $A$, there exists a finite sequence
$i=i_0, i_1, \cdots, i_n=j$ such that
$i\in \bar{J}\cup\{i\;|\; a_{i,i}=0\}$,
$i_k \not\in \bar{J}\cup\{i\;|\; a_{i,i}=0\}$ for $k\geq1$,
and $a_{i_k, i_{k+1}}\neq0$ for all $k$.
What we have found above allows us to recursively show
$\im(\xi)\subset T_{i_k}$ and in particular,
$\im(\xi)\subset T_j$.
\end{proof}

\begin{prop}\label{6.7}
Suppose there exists some finite $J\subset\I$ such that for every $j\in J$,
$a_{j,j} = 2$ and for which the corresponding submatrix
$A_J = (a_{i,j})_{i,j\in J}$ is
indecomposable and not of finite type. Then, $\cntr_J\subset U^0$.
\end{prop}
\begin{proof} By \cite[Proposition 4.9]{Kac}, we have $|W_J|=\infty$.
Let $J'\subset J$ be such that $|W_{J'}|=\infty$ and $|W_{J''}|\lneq\infty$
for all $J''\subsetneq J'$.
We may use Proposition~\ref{mulguisin}
if we can show $\cntr_{J'}\subset U^0$. Hence it suffices to show that
if $h\in P^\vee$, $|W_{J'}(h)| \lneq\infty$, then $\ali(h)=0$,
for all $i\in J'$.
Give partial order to $\h$ by setting $h_1\geq h_2$ if and only if
$h_1-h_2\in (\sum_i \Z_{\geq0}h_i) + (\sum_i \Z_{\geq0}d_i)$.
Let $h'\in W_{J'}(h)$ be maximal with respect to this order.
Then for each $i\in J'$, if $\ali(h')<0$, then $h'< r_i h'$, so
$\ali(h')\geq0$ for all $i\in J'$.
Set $W_{h'} = \{ w\in W\;|\; w(h') = h' \}$.
By \cite[Proposition 3.12(a)]{Kac}, $W_{h'} = W_{J''}$ with
$J'' = \{ i\in J \;|\; \ali(h')=0\}$. If $J''\subsetneq J'$, then
$|W_{J'}(h')| = | W_{J'} / W_{h'} | = \infty$. Hence we must have
$J'' = J'$ and $\ali(h')=0$ for all $i\in J'$.
$\{h'\} = W_{J'}(h') = W_{J'}(h)$. So $h=h'$ and $\ali(h) = 0$ for
all $i\in J'$.
\end{proof}

\begin{prop}\label{6.8}
Let $A$ be indecomposable, not of finite type, and
$a_{i,i} = 2$ for all $i\in\I$.
Then, $\cntr\subset U^0$.
\end{prop}
\begin{proof}
Suppose there exists some finite indecomposable submatrix which is not
of finite type. Then we may use Proposition~\ref{mulguisin} and
Proposition~\ref{6.7} to obtain the result.

If, to the contrary, every finite submatrix of $A$ is of finite type, 
it must be one of the following types :
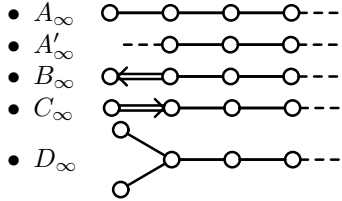
\begin{itemize}
\item $A_\infty$\quad
      \begin{texdraw}
      \rcir\rlvec(6 0)\rcir\rlvec(6 0)\rcir\rlvec(6 0)\rcir
      \lpatt(1 1)\rlvec(6 0)
      \end{texdraw}
\item $A_\infty'$\quad
      \begin{texdraw}
      \move(8 0)\lpatt(1 1)\rlvec(-6 0)\move(8 0)\lpatt()
      \rcir\rlvec(6 0)\rcir\rlvec(6 0)\rcir
      \lpatt(1 1)\rlvec(6 0)
      \end{texdraw}
\item $B_\infty$\quad
      \begin{texdraw}
      \rcir\move(2 0.3)\rlvec(6 0)\move(8 -0.3)\rlvec(-6 0)
      \move(2 0)\rlvec(1.2 1.2)\move(2 0)\rlvec(1.2 -1.2)\move(8 0)
      \rcir\rlvec(6 0)\rcir\rlvec(6 0)\rcir
      \lpatt(1 1)\rlvec(6 0)
      \end{texdraw}
\item $C_\infty$\quad
      \begin{texdraw}
      \rcir\move(2 0.3)\rlvec(6 0)\move(8 -0.3)\rlvec(-6 0)
      \move(8 0)\rlvec(-1.2 1.2)\move(8 0)\rlvec(-1.2 -1.2)\move(8 0)
      \rcir\rlvec(6 0)\rcir\rlvec(6 0)\rcir
      \lpatt(1 1)\rlvec(6 0)\move(-0.4 0)
      \end{texdraw}
\item $D_\infty$\quad
      \raisebox{-0.4\height}
      {\begin{texdraw}
      \move(8 0)\rcir\rlvec(6 0)\rcir\rlvec(6 0)\rcir
      \lpatt(1 1)\rlvec(6 0)\lpatt()
      \move(9 0)\rmove(-0.86 0.5)\rlvec(-5.1 3)\rmove(-0.86 0.5)\lcir r:1
      \move(9 0)\rmove(-0.86 -0.5)\rlvec(-5.1 -3)\rmove(-0.86 -0.5)\lcir r:1
      \end{texdraw}}
\end{itemize}
In all cases, with $\I$ naturally ordered, the matrix satisfies the
following condition.
\begin{quote}\label{property}
For each $i\in\I$, there exists some $j > i$
such that $a_{i,j}\neq0$, and $a_{i,k}=0$ for $k>j$.
\end{quote}
Let $c = \sum_h c_h q^h\in \im(\xi) \subset (U^0)^{\widetilde{W}}$.
Fix $h\in P^\vee$ for which $c_h\neq0$. We aim to show
$|\widetilde{W}(h)| = \infty$ if $\alj(h)\neq0$ for some $j\in\I$.
We may assume that only finitely many $j\in\I$ satisfy $\alj(h)\neq0$.
Let $k\in\I$ be the maximal of those so that $\alj(h)=0$ for all $j>k$
and $\alpha_k(h)\neq0$.
Set $i_0 = k$, and using property~\eqref{property}, recursively choose $i_n$
so that $i_{n+1} > i_n$ and $a_{{i_n},{i_{n+1}}}\neq0$. Put $h_0 = h$ and
$h_{n+1} = r_{i_n} h_n$. Then, $h_n$ cannot form a closed orbit and
$|\widetilde{W}(h)| = \infty$. Hence
$\im(\xi)\subset \oplus_{h\in P^\vee, \ali(h)=0} \C q^h$
and $\cntr\subset U^0$.
\end{proof}

We can now collect all results and state :
\begin{thm}
Assume that the Borcherds-Cartan matrix $A = (a_{i,j})_{i,j\in\I}$ is
indecomposable and not of finite type.
Suppose that every $0_j^-$ is connected directly to a
$0_i$ or a $\circleddash_i$. Except for the case when
$|\I|=1$ with $m_i=1$, the center $\cntr$ belongs to $U^0$.
\end{thm}
\begin{proof}
We apply Proposition~\ref{mulguisin} to each possible case.

If $|\I|=1$, the conditions imply either a $0_i^+$ or a
$\circleddash_i$ with $m_i \geq 2$. These cases may be handled by
Proposition~\ref{6.2}~\ref{gau} and Lemma~\ref{rank2}~\ref{ra}, respectively.

Now suppose $|\I|\geq2$. Proposition~\ref{6.8} does away with the case when
all $a_{i,i}=2$. If it contains a $0_i^+$, we may again use
Proposition~\ref{6.2}~\ref{gau}.
If it contains a $0_i^-$ but no $0_i^+$, we use
Lemma~\ref{rank2}~\ref{rf},\ref{re}. The only other case
is covered by Lemma~\ref{rank2}~\ref{rc},\ref{rb}.
\end{proof}
%

\section{The Universal R-matrix}

In this section, we find the universal R-matrix for the quantum group
$U_q(\g)$.

A Hopf superalgebra (a colored Hopf algebra)
$\Hf$ together with an element
$\R\in \Hf\otimes\Hf$ is called a
\defi{quasi-triangular Hopf superalgebra} if it satisfies:
\begin{enumerate}
\item $\R$ is invertible,
\item $\R\cdot\Delta(a) = \Delta'(a)\cdot\R$ for all $a\in\Hf$,
\item $(\Delta\otimes1)(\R) = \R_{13}\R_{23}$,
\item $(1\otimes\Delta)(\R) = \R_{13}\R_{12}$,
\end{enumerate}
where $\Delta' = P\circ\Delta$ with $P$ a colored twisting map,
and where $\R_{ij}$ is an element of
$\Hf\otimes\Hf\otimes\Hf$ such that the $i$'th and $j$'th components are
given by $\R$ and the remaining component is 1.
The element $\R$ is called the \defi{universal R-matrix}.
It satisfies the Yang-Baxter equation
\begin{equation}\label{YBE}
\R_{12}\R_{13}\R_{23} = \R_{23}\R_{13}\R_{12}.
\end{equation}

A Hopf superalgebra $\Hf$ together with an element $\Ct\in\Hf\otimes\Hf$
and an algebra homomorphism
$\Phi : \Hf\otimes\Hf \longrightarrow \Hf\otimes\Hf$ is called a
\defi{pre-triangular Hopf superalgebra} if it satisfies:
{ 
\renewcommand{\theenumi}{\arabic{enumi}}
\renewcommand{\labelenumi}{\textnormal{(P\theenumi)}}
\begin{enumerate}
\item $\Ct$ is invertible,
\item $\Ct\cdot\Delta(a) = \Phi(\Delta'(a))\cdot\Ct$ for all $a\in\Hf$,
\item $\Phi_{23}\circ\Phi_{13}(\Ct_{12}) = \Ct_{12}$,
\item $\Phi_{12}\circ\Phi_{13}(\Ct_{23}) = \Ct_{23}$,
\item $\Phi_{23}(\Ct_{13})\cdot\Ct_{23} = (\Delta\otimes1)(\Ct)$,
\item $\Phi_{12}(\Ct_{13})\cdot\Ct_{12} = (1\otimes\Delta)(\Ct)$.
\end{enumerate}
} 
Under some conditions, it is possible to show that a pre-triangular
Hopf superalgebra becomes a quasi-triangular Hopf superalgebra.

We set $U^{+,\beta} = %
\bigoplus_{\gamma\in Q^+, \gamma\nleq\beta}U^+_\gamma$ for each $\beta\in Q^+$
and define the completion $\widehat{U}$ of $U$ by:
\begin{equation}
\widehat{U} = \lim_{\underset{\beta}{\longleftarrow}}\, U/UU^{+,\beta}.
\end{equation}
There is a natural embedding of $U$ in $\widehat{U}$ and there is a natural
algebra structure on $\widehat{U}$ which extends that of $U$ under this
embedding.

The completion of $U^{\otimes n}$ is similarly defined. We will write
$\widehat{U}\hat{\otimes}\widehat{U}$ for the completion of $U\otimes U$.

Define an algebra automorphism
$\Phi : U\otimes U \longrightarrow U\otimes U$ by
\begin{alignat}{2}
\Phi(q^h\otimes q^{h'}) &= q^h\otimes q^{h'},\\
\Phi(\eik\otimes1) &= \eik\otimes K_i,&\quad
\Phi(1\otimes\eik) &= K_i\otimes\eik,\\
\Phi(\fik\otimes1) &= \fik\otimes K_i^{-1},&
\Phi(1\otimes\fik) &= K_i^{-1}\otimes\fik.
\end{alignat}
It can be shown that $\Phi$ naturally extends to an algebra automorphism of
$\widehat{U}\hat{\otimes}\widehat{U}$.

We denote by $\Ct_\beta\in U^+_\beta\otimes U^-_{-\beta}$ the canonical
element of the bilinear form
$(\ |\ ) : U^+_\beta\times U^-_{-\beta} \longrightarrow \C$.
Define
\begin{equation}
\Ct = \sum_{\beta\in Q^+} \theta(\beta,\beta)q^{(h_\beta|h_\beta)}
(K_\beta^{-1}\otimes K_\beta)\Ct_\beta\in \widehat{U}\hat{\otimes}\widehat{U}.
\end{equation}

\begin{lem}\hfill
\begin{enumerate}
\item $\Ct\cdot\Delta(q^h) = \Phi(\Delta'(q^h))\cdot\Ct \quad(h\in P^\vee)$,
\item $(\Phi_{23}\circ\Phi_{13})(\Ct_{12}) = \Ct_{12}$,
\item $(\Phi_{12}\circ\Phi_{13})(\Ct_{23}) = \Ct_{23}$.
\end{enumerate}
\end{lem}
\begin{proof}
This is just straightforward calculation.
\end{proof}

\begin{lem}\label{lm:72}
Let $\beta\in Q^+$.
\begin{enumerate}
\item $\sum_{\substack{\gamma,\delta\in Q^+\\ \gamma+\delta = \beta}}
  \Ct_\gamma(K_\delta\otimes1)(S\otimes1)(\Ct_\delta) = \delta_{\beta,0}$.
\label{lema}
\item $\sum_{\substack{\gamma,\delta\in Q^+\\ \gamma+\delta = \beta}}
  (K_\gamma\otimes1)(S\otimes1)(\Ct_\gamma)\Ct_\delta = \delta_{\beta,0}$.
\item $\theta_{i,i}[1\otimes\eik,\Ct_{\beta+\ali}] =
  C_\beta(\eik\otimes K_i^{-1}) - (\eik\otimes K_i)\Ct_\beta$.
\item $\theta_{i,i}[\fik\otimes1,\Ct_{\beta+\ali}] =
  C_\beta(K_i\otimes\fik) - (K_i^{-1}\otimes\fik)\Ct_\beta$.
\item $(\Delta\otimes1)(\Ct_\beta) =
  \sum_{\substack{\gamma,\delta\in Q^+\\ \gamma+\delta = \beta}}
  q^{-(h_\gamma|h_\delta)}(K_\delta\otimes1\otimes1)
  (\Ct_\gamma)_{13}(\Ct_\delta)_{23}$.
\item $(1\otimes\Delta)(\Ct_\beta) =
  \sum_{\substack{\gamma,\delta\in Q^+\\ \gamma+\delta = \beta}}
  q^{-(h_\gamma|h_\delta)}(1\otimes1\otimes K_{-\delta})
  (\Ct_\gamma)_{13}(\Ct_\delta)_{12}$.
\end{enumerate}
\end{lem}
\begin{proof}
Here we show the proof for \ref{lema} only.
Other cases may be proved in a similar spirit.

The case $\beta=0$ is trivial.
So assume $\beta\in Q^+ \setminus\{0\}$.
The left hand side is contained in $U^+_\beta\otimes U$, so by
Theorem~\ref{nondg} it suffices to
show that the application of $(\cdot|w)\otimes1$ is zero for all
$w\in U^-_{-\beta}$.
We may write
\begin{equation*}
\Delta(w) = \sum_{\substack{\gamma,\delta\in Q^+ \\ \gamma+\delta = \beta}}
w^{\delta,\gamma}(1\otimes K_\delta)
\quad\text{with}\quad
w^{\delta,\gamma} \in U^-_{-\delta}\otimes U^-_{-\gamma}
\end{equation*}
and
\begin{equation*}
w^{\delta,\gamma} = \sum_m
w^{\delta,\gamma}_{\delta,m}\otimes w^{\delta,\gamma}_{\gamma,m}
\quad\text{with}\quad
w^{\delta,\gamma}_{\delta,m}\in U^-_{-\delta},\quad
w^{\delta,\gamma}_{\gamma,m}\in U^-_{-\gamma}.
\end{equation*}
We may also fix basis $\{x^\gamma_r\}_r$ and $\{y^\gamma_r\}_r$ of 
$U^+_\gamma$ and $U^-_{-\gamma}$, respectively, which are dual
with respect to the bilinear form.
Now,
{\allowdisplaybreaks
\begin{align*}
((\ |w)\otimes&1)(\text{LHS}) =
((\ |w)\otimes1)\Big(
\sum_{\gamma,\delta,r,s}
\theta(y^\gamma_r,x^\delta_s)x^\gamma_r K_\delta S(x^\delta_s)
\otimes
y^\gamma_r y^\delta_s
\Big) \\
&=\sum_{\gamma,\delta,r,s}
\theta(x^\gamma_r, K_\delta S(x^\delta_s))\theta(-\gamma,\delta)
(K_\delta S(x^\delta_s) \otimes x^\gamma_r | \Delta(w))
y^\gamma_r y^\delta_s\\
&=\sum_{\gamma,\delta,r,s,m}
\theta(\delta,\gamma)
(k_\delta S(x^\delta_s)|w^{\delta,\gamma}_{\delta,m})
(x^\gamma_r|w^{\delta,\gamma}_{\gamma,m}K_{-\delta})
y^\gamma_r y^\delta_s\\
&= \sum_{\gamma,\delta,m}\theta(\delta,\gamma)
\Big( \sum_r
(x^\gamma_r|w^{\delta,\gamma}_{\gamma,m}K_{-\delta})y^\gamma_r
\Big)\Big( \sum_s
(K_\delta S(x^\delta_s)|w^{\delta,\gamma}_{\delta,m})y^\delta_s
\Big)\\
&= \sum_{\gamma,\delta,m}\theta(\delta,\gamma)
\Big( \sum_r
(x^\gamma_r|w^{\delta,\gamma}_{\gamma,m})y^\gamma_r
\Big)\Big( \sum_s
((x^\delta_s|K_\delta^{-1} S^{-1}(w^{\delta,\gamma}_{\delta,m}))y^\delta_s
\Big)\\
&= \sum_{\gamma,\delta,m}\theta(\delta,\gamma)
w^{\delta,\gamma}_{\gamma,m}K_\delta^{-1}S^{-1}(w^{\delta,\gamma}_{\delta,m})\\
&= S^{-1}\Big(
\sum_{\gamma,\delta,m}
w^{\delta,\gamma}_{\delta,m}
S(w^{\delta,\gamma}_{\gamma,m}K_\delta^{-1})
\Big)\\
&= (S^{-1}\circ \text{m} \circ (1\otimes S)\circ\Delta)(w)\\
&= \varepsilon(w)\\
&= 0.
\end{align*}
} 
Hence the left hand side is zero when $\beta \neq 0$.
\end{proof}

\begin{prop}
Let
\begin{equation}
\Ct' = \sum_{\beta\in Q^+} \theta(\beta,\beta)q^{(h_\beta|h_\beta)}
(1\otimes K_\beta)(S\otimes1)(\Ct_\beta)\in\widehat{U}\hat{\otimes}\widehat{U}.
\end{equation}
Then $\Ct\Ct'=\Ct'\Ct=1$.
\end{prop}
\begin{proof}
{\allowdisplaybreaks
\begin{align*}
\Ct\Ct' &=
\Big(\sum_{\gamma\in Q^+}\theta(\gamma,\gamma)q^{(h_\gamma|h_\gamma)}
(K_\gamma^{-1}\otimes K_\gamma)\Ct_\gamma
\Big)\\
&\phantom{\sum_{\gamma\in Q^+}\theta(\gamma,\gamma)q^{(h_\gamma|h_\gamma)}}
\cdot\Big(
\sum_{\delta\in Q^+}\theta(\delta,\delta)q^{(h_\delta|h_\delta)}
(1\otimes K_\delta)(S\otimes1)(\Ct_\delta)
\Big)\\
&= \sum_{\beta\in Q^+}
\sum_{\substack{\gamma+\delta=\beta\\ \gamma,\beta\in Q^+}}
\theta(\beta,\beta)q^{(h_\gamma|h_\gamma)+(h_\delta|h_\delta)}
(K_\gamma^{-1}\otimes K_\gamma)\Ct_\gamma(1\otimes K_\delta)
(S\otimes1)(\Ct_\delta)\\
&= \sum_{\beta\in Q^+}
\theta(\beta,\beta)q^{(h_\beta|h_\beta)}
(K_\beta^{-1}\otimes K_\beta)
\sum_{\substack{\gamma+\delta=\beta\\ \gamma,\beta\in Q^+}}
\Ct_\gamma(K_\delta\otimes1)(S\otimes1)(\Ct_\delta)
\end{align*}
}
We may now apply Lemma~\ref{lm:72}.
The other part is done similarly.
\end{proof}

\begin{prop}
We have
\begin{align}
\Ct\cdot\Delta(\eik) &= \Phi(\Delta'(\eik))\cdot\Ct,\\
\Ct\cdot\Delta(\fik) &= \Phi(\Delta'(\fik))\cdot\Ct.
\end{align}
\end{prop}
\begin{proof}
\begin{align*}
\Ct\cdot\Delta(\eik) &=
\sum_{\beta\in Q^+}
\theta(\beta,\beta)q^{(h_{\beta+\ali}|h_\beta)}
(K_\beta^{-1}\otimes K_{\beta+\ali})\Ct_\beta
(\eik\otimes K_i^{-1})\\
&\phantom{\sum_{\beta\in Q^+}\theta(\beta,\beta)q} +
\sum_{\beta\in Q^+}
\theta(\beta,\beta)q^{(h_{\beta-\ali}|h_\beta)}
(K_{\beta-\ali}^{-1}\otimes K_\beta)\Ct_\beta
(1\otimes\eik)
\end{align*}
\begin{align*}
\Phi(\Delta'(\eik))\cdot\Ct =
\sum_{\beta\in Q^+}
&\theta(\beta,\beta)q^{(h_\beta|h_{\beta-\ali})}
(K_{\beta-\ali}^{-1}\otimes K_\beta)(1\otimes\eik)\Ct_\beta\\
&+ \sum_{\beta\in Q^+}
\theta(\beta,\beta)q^{(h_{\beta+\ali}|h_\beta)}
(K_\beta^{-1}\otimes K_{\beta+\ali})(\eik\otimes K_i)\Ct_\beta
\end{align*}
\begin{align*}
\Ct\cdot&\Delta(\eik) - \Phi(\Delta'(\eik))\cdot\Ct\\
&= \sum_{\beta\in Q^+}
\theta(\beta,\beta)q^{(h_\beta|h_{\beta+\ali})}
(K_\beta^{-1}\otimes K_{\beta+\ali})\\
&\phantom{\sum\theta(\gamma,\gamma)}\cdot\Big\{
\Ct_\beta(\eik\otimes K_i^{-1})-(\eik\otimes K_i)\Ct_\beta
- \theta(\ali,\ali)[1\otimes\eik, \Ct_{\beta+\ali}]
\Big\}
\end{align*}
We apply Lemma~\ref{lm:72} to obtain the result.
The other case is similar.
\end{proof}

\begin{prop}
We have
\begin{align}
\Phi_{23}(\Ct_{13})\cdot\Ct_{23} &= (\Delta\otimes1)(\Ct),\\
\Phi_{12}(\Ct_{13})\cdot\Ct_{12} &= (1\otimes\Delta)(\Ct).
\end{align}
\end{prop}
\begin{proof}
\begin{equation*}
\Phi_{23}(\Ct_{13}) =
\sum_{\gamma\in Q^+}
\theta(\gamma,\gamma) q^{(h_\gamma|h_\gamma)}
(K_\gamma^{-1}\otimes K_\gamma^{-1}\otimes K_\gamma)(\Ct_\gamma)_{13}
\end{equation*}
\begin{align*}
(\Phi_{23}(\Ct_{13}))\Ct_{23} =
\sum_{\beta\in Q^+}
\theta(\beta,\beta)&q^{(h_\beta|h_\beta)}
(K_\beta^{-1}\otimes K_\beta^{-1}\otimes K_\beta)\\
&\cdot\sum_{\substack{\gamma,\delta\in Q^+\\ \gamma+\delta = \beta}}
q^{-(h_\delta|h_\gamma)}(K_\delta\otimes1\otimes1)
(\Ct_\gamma)_{13}(\Ct_\delta)_{23}
\end{align*}
\begin{equation*}
(\Delta\otimes1)(\Ct) =
\sum_{\beta\in Q^+}
\theta(\beta,\beta)q^{(h_\beta|h_\beta)}
(K_\beta^{-1}\otimes K_\beta^{-1}\otimes K_\beta)
(\Delta\otimes1)(\Ct_\beta)
\end{equation*}
The second case is done similarly.
\end{proof}

The Propositions tell us that $U$ is almost a pre-triangular Hopf
superalgebra.
\begin{thm}
The statements \textnormal{(P1)} and \textnormal{(P2)} hold in
$\widehat{U}\hat{\otimes}\widehat{U}$ and
the relations \textnormal{(P3)}--\textnormal{(P6)} hold in
$\widehat{U}\hat{\otimes}\widehat{U}\hat{\otimes}\widehat{U}$.
\end{thm}

A weight module is \defi{$P$-weighted} if all its weights belong to $P$.
Notice $(P|P) \subset \Z$.
This allows us to define $\mathfrak{Z}\in\text{End}(V\otimes W)$
for any $P$-weighted $U_q(\g)$-modules $V$ and $W$ by setting,
\begin{equation}
\mathfrak{Z}( v\otimes w ) = 
q^{(\text{wt}(v)|\text{wt}(w))} v\otimes w
\end{equation}
on homogeneous elements and extending by linearity.
The map $\mathfrak{Z}$ is certainly invertible.
There is a natural action of $U\otimes U$
on $V\otimes W$ and as endomorphisms on
$V\otimes W$,
\begin{equation}\label{Phi-Z}
\Phi(a\otimes b) = \mathfrak{Z}\circ(a\otimes b)\circ\mathfrak{Z}^{-1}
\end{equation}
for every $a\otimes b\in U\otimes U$.

Set $\R = \mathfrak{Z}^{-1}\Ct$. Then we finally have:
\begin{thm}
Let $V_i$ \textup{(}$i=1,2,3$\textup{)} be $P$-weighted $U_q(\g)$-modules.
As endomorphisms on $V_1\otimes V_2\otimes V_3$,
when it can be defined, $\R$ satisfies the Yang-Baxter equation~\eqref{YBE}.
\end{thm}
\begin{proof}
 From (P5) and equation~\eqref{Phi-Z}, we have
\begin{gather}
\ZP_{23}\Ct_{13}\ZP_{23}^{-1}\Ct_{23} = (\Delta\otimes1)(\Ct)\\
\R_{13}\R_{23} = \ZP_{13}^{-1}\ZP_{23}^{-1}(\Delta\otimes1)(\Ct)\label{eqtmp}
\end{gather}
Applying $P\otimes1$ to both sides of (P5) and working as above, we get
\begin{equation}
\R_{23}\R_{13} = \ZP_{23}^{-1}\ZP_{13}^{-1}(\Delta'\otimes1)(\Ct).
\end{equation}
The use of (P2) shows,
\begin{align}
\R_{23}\R_{13}\R_{12} &=
     \ZP_{23}^{-1}\ZP_{13}^{-1}((\Delta'\otimes1)(\Ct))\R_{12}\\
  &= \ZP_{23}^{-1}\ZP_{13}^{-1}%
			((\Delta'\otimes1)(\Ct))(\ZP^{-1}\Ct\otimes1)\\
  &= \ZP_{23}^{-1}\ZP_{13}^{-1}(\ZP^{-1}\Ct\otimes1)(\Delta\otimes1)(\Ct)
\end{align}
Now, the $\ZP_{ij}$ commute with each other and (P3) with \eqref{Phi-Z}
says $\Ct_{12}$ commutes with $\ZP_{13}^{-1}\ZP_{23}^{-1}$,
so we may use \eqref{eqtmp} to write
\begin{equation}
\ZP_{23}^{-1}\ZP_{13}^{-1}(\ZP^{-1}\Ct\otimes1)(\Delta\otimes1)(\Ct) =
\R_{12}\R_{13}\R_{23}.
\end{equation}
Putting things together, we have the result.
\end{proof}

\emph{Acknowledgments} I would like to thank Professor Seok-Jin Kang and one
of my colleagues, Kyu-Hwan Lee, for help throughout this work.
Lee has contributed much in finding the structure theorem for the center
of finite type cases.

\providecommand{\bysame}{\leavevmode\hbox to3em{\hrulefill}\thinspace}

\end{document}